\documentclass[12pt,twoside,a4paper,notitlepage]{article}
\usepackage[latin1]{inputenc}
\usepackage{color}
\usepackage{amsmath}
\usepackage{amsfonts}
\usepackage{amssymb}
\usepackage{graphicx}
\usepackage[normalem]{ulem}

\usepackage[a4paper,margin=2cm,vmargin={1cm,2cm},includeheadfoot]{geometry}

\newtheorem{theorem}{Theorem}[section]

\newcommand{\qed}{\hspace{\stretch{3}}$\square$\vspace{1.8ex}}

\newtheorem{lemma}[theorem]{Lemma}

\newtheorem{prop}[theorem]{Proposition}
\newtheorem{teo}[theorem]{Theorem}
\newtheorem{remark}[theorem]{Remark}

\def\E{\mathbb{E}}

\def\N{\mathbb{N}}

\def\V{\mathbb{V}}

\newcommand{\proof}{{\bf Proof. }}
\newcommand{\ie}{i.e. }

\title{\bf  Multilevel Monte Carlo simulation for L\'evy processes based on the Wiener-Hopf factorisation}
\author{\small{\sc A. Ferreiro-Castilla$^{\dagger,}$\footnote{Supported by a Royal Society Newton International Fellowship.}  \ \ A.E. Kyprianou\footnote{Department of Mathematical Sciences, University of Bath,
Claverton Down, Bath, BA2 7AY, U.K.} 
\ \ R. Scheichl$^\dagger$%l\footnote{Department of Mathematical Sciences, University of Bath, Claverton Down, Bath, BA2 7AY, U.K.} 
\ \
  G. Suryanarayana\footnote{Numerical Analysis and Applied Mathematics Section, KU Leuven,
Celestijnenlaan 200a,  box 2402,
3001 Heverlee, Leuven, Belgium.}
}
}

\begin{document}

\maketitle
\begin{abstract}
In Kuznetsov et al. \cite{KKPvS} a new Monte Carlo simulation technique was introduced for a large family of L\'evy processes that is based on the Wiener-Hopf decomposition. We pursue this idea further by combining their technique with the recently introduced multilevel Monte Carlo methodology. Moreover, we provide here for the first time a theoretical analysis of the new Monte Carlo simulation technique in \cite{KKPvS} and of its multilevel variant for computing expectations of functions depending on the historical trajectory of a L{\'e}vy process. We derive rates of convergence for both methods and show that they are uniform with respect to the ``jump activity'' (e.g. characterised by the Blumenthal-Getoor index). We also present a modified version of the algorithm in Kuznetsov et al. \cite{KKPvS} which combined with the multilevel methodology obtains the optimal rate of convergence for general L\'evy processes and Lipschitz functionals. This final result is only a theoretical one at present, since it requires independent sampling from a triple of distributions which is currently only possible for a limited number of processes.

\bigskip

\bigskip


\end{abstract}

\section{Introduction}

A classical problem in mathematical finance deals with the ability to compute $\E[F(X)]$, where $X:=\{X_{s}: s\in [0,t]\}$ is a stochastic process which models the underlying risky asset and $F$ is a payoff function which may depend on the historical path of $X$.  This is the typical setting for pricing a variety of exotic options in finance such as look back and barrier options.  Starting with the early work of Madan and Seneta \cite{MD}, a class of processes playing the role of $X$ that have found prominence in this respect is the class of L\'evy processes. An extensive overview of their position in mathematical finance can be found in the books \cite{BL, CT, S, SC}. More recently L\'evy processes have also been extensively used in modern insurance risk theory; see for example Asmussen and Albrecher \cite{AA} and Kl\"uppelberg et al. \cite{KKM}. In insurance mathematics, it is  the L\'evy process itself which models the surplus wealth of an insurance company until its ruin. There are also extensive applications of L\'evy processes in queuing theory, genetics and mathematical biology as well as in stochastic differential equations {(see e.g.~\cite{CLU,D11,Doe,KW})}. 

In both financial and insurance settings, a  key quantity of generic interest is the joint law of the current position and the running maximum of a L\'evy process at a fixed time, if not the individual marginals associated with the latter bivariate law. Consider the following example. If we define $\overline{X}_t = \sup_{s\leq t}X_s$ then the pricing of barrier options boils down to evaluating expectations of the form $\mathbb{E}[f(x+X_t) \mathbf{1}_{\{x+\overline{X}_t >b\}}]$ for some threshold $b>0$. Indeed if $f(x) = (K-{\mathrm e}^x)^+$ and $K>0$, then the latter expectation is related to the value of an ``up-and-in" put. In credit risk one is predominantly interested  in the quantity $\widehat{\mathbb{P}}(\overline{X}_t <x)$ as a function in $x$ and~$t$, where $\widehat{\mathbb{P}}$ is the law of the dual process $-X$. Indeed  it is a functional of the latter probabilities that gives the price of a credit default swap {not to mention the recently introduced financial instruments known as convertible contingencies (CoCos). See for example the recent book of Schoutens and Cariboni \cite{SC} as well as Corcuera et al. \cite{coco}}. One is similarly interested in $\widehat{\mathbb{P}}(\overline{X}_t \geq x)$ in ruin theory, since these probabilities are also equivalent to the finite-time ruin probabilities; {cf. Asmussen and Albrecher \cite{AA}}. 

A widely used approach to compute expectations of functions depending on the historical trajectory of a L{\'e}vy process over the time horizon, say $[0,t]$, is to approximate the path by a random walk with $n$ steps, each step covering $t/n$ units of time, and therewith to perform a Monte Carlo (MC) simulation. Giles \cite{Giles07,Giles08} introduced an adaptation of the straightforward MC methodology, the multilevel Monte Carlo method (MLMC), which is especially suited to the scenario we are interested in, but  in the case that  $X$ is a pure diffusion process. Very recently, there has been increasing attention to the MLMC method also in the setting of L{\'e}vy processes, see Giles and Xia \cite{GX11} for the jump diffusion setting, and Dereich \cite{D11} and Dereich and  Heidenreich \cite{DH11} for more general L\'evy processes. Generally speaking all these methods share a common approach, which consists in constructing an embedded sequence of grids that are made up of a mixture of deterministic and random points. The random points in these grids deal with the large jumps of the L\'evy process and the deterministic points deal with the ``small movements", that is to say, the diffusive part and/or the small jumps.

In this paper we consider an alternative method based entirely on a random grid. In particular we shall introduce an adaptation of the MLMC method based on a very recently introduced technique for performing MC simulations that appeals to the so-called Wiener-Hopf factorisation for one-dimensional L{\'e}vy processes. This last technique is called the Wiener-Hopf Monte Carlo (WHMC) simulation method and it was introduced in Kuznetsov et al. \cite{KKPvS}. We will denote the coupling of these two techniques the multilevel Wiener-Hopf Monte Carlo method (MLWH). Our analysis will focus on the particular setting that $X$ is a one-dimensional L{\'e}vy process and $F$ is a Lipschitz function which depends  on the  value of $X$ and its past supremum at a fixed time $t>0$. 

The main result in this paper shows that the root mean square error in the MLWH method converges in general with order $\mathcal{O}(\nu^{-\frac{1}{4}} )$ for processes of unbounded variation and $\mathcal{O}(\nu^{-\frac{1}{3}} )$ for processes of bounded variation. Unlike Dereich~\cite{D11}, Dereich and Heidenreich \cite{DH11}, and Giles and Xia \cite{GX11}, our method is robust in the sense that the convergence rate does not exhibit dependence on the nature of the jump structure (in particular small jumps) of the underlying L\'evy process. This means that we can ensure a convergence rate of at least $\mathcal{O}(\nu^{-\frac{1}{4}} )$ independently of the structure of the L\'evy measure, which is important in the implementation of numerical methods for calibration purposes, as slightly different parameters in a model can lead to completely different behavior of the process with respect to the jump structure. Moreover, the implementation of the algorithm is extremely straightforward, requiring only to be able to independently sample from two distributions. This robustness comes at the price of our algorithm not having  optimal convergence rates, in particular when the underlying L\'evy process has paths of bounded variation.  For L\'evy processes of unbounded variation we are able to  derive regimes where our methodology performs better than the algorithms in Dereich~\cite{D11}, Dereich and Heidenreich \cite{DH11} and Jacod et al. \cite{JKMP05} and, as far as we know, better than  any other existing algorithm. We will also show the existence of another multi-level algorithm, which is a modification of the original algorithm in Kuznetsov et al. \cite{KKPvS}, that achieves a convergence order $\mathcal{O}(\nu^{-\frac{1}{2}} \log^2 \nu)$ independently of the jump structure. As shown in Creutzig et al.~\cite{CDMR09}, up to the logarithmic term this is in fact optimal for general L\'evy processes. Unfortunately however, this algorithm requires sampling from a triplet of dependent distributions for which there are currently very few known examples of L\'evy processes where this is possible. We shall elaborate further on this point later in this paper.

As alluded to above, we shall benchmark our  MLWH results against those of  Dereich and Heidenreich \cite{DH11} and Dereich \cite{D11}. Note that, even though these two papers consider the more general setting of expectations of path functionals of solutions to  SDEs driven by general L\'evy processes, the scenario we consider here is permitted within their framework and we are currently working on extending our methodology to the more general setting. 

The paper is organised as follows. In the next section we will review the general setting for the Wiener-Hopf factorisation of L{\'e}vy processes and describe the  WHMC method introduced in Kuznetsov et al. \cite{KKPvS}. Thereafter, in Section \ref{MLWH_sec}, we describe how an adaptation of the WHMC method can be used in the context of multilevel simulation, that is to say, we formally introduce the MLWH method. Section \ref{num-anal} is devoted to the numerical analysis of the one level (WHMC) and the multilevel (MLWH) method, in particular, giving mathematical justification to our claims concerning the performance of the MLWH method above. Moreover, we discuss how the MLWH method performs against other approaches that have been mentioned. Section \ref{randomisation_time} will introduce a theoretical variation of the original MLWH, which we can show achieves optimal convergence rates but which, unfortunately, cannot  yet be implemented in practice until more substantial examples  of the Wiener-Hopf factorisation are discovered. Finally, Section \ref{Numerics} provides some implementation details and numerical computations for a representative example which support the previous theoretical derivations and confirm the feasibility of the approach.

\section{Preliminaries}\label{WHMC_MLMC}

\subsection{L\'evy processes and the Wiener-Hopf factorisation}

We begin  by briefly reviewing the definition of a one-dimensional L\'evy process and its associated Wiener-Hopf factorisation. For a more in-depth account we refer the reader to the monographs of Bertoin \cite{bertoin}, Kyprianou \cite{kyp-book} or Sato \cite{Sato99}. Note that the Wiener-Hopf factorisation only exists for one-dimensional L\'evy processes.

Recall that a one-dimensional L\'evy process with law $\mathbb{P}$, henceforth denoted by $X := \{X_t : t\geq 0\}$, is a stochastic process issued from the origin which enjoys the properties of having stationary and independent increments with paths that are almost surely right-continuous with left limits. It is a well understood fact that, as a consequence of  this definition,  the law of every L\'{e}vy process is characterised through a triplet $(a,\sigma,\Pi)$, where $a\in\mathbb{R}$, $\sigma\geq0$ and $\Pi$ is a measure concentrated on $\mathbb{R}\backslash\{0\}$ such that $\int_{\mathbb{R}}(1\wedge x^{2})\Pi(\mathrm{d}x)<\infty$. More precisely, for all $t \ge 0$, 
\[
\mathbb{E}[{\mathrm e}^{{\rm i}\theta X_{t}}]={\mathrm e}^{-t\Psi(\theta)}, \qquad \text{for all} \ \ \theta\in\mathbb{R},
\]
where
\begin{equation*}
\Psi ( \theta )  \ =\  \mathrm{i}a\theta +\frac{1}{2}\sigma
^{2}\theta
^{2}+\int_{\mathbb{R}}(1-\mathrm{e}^{\mathrm{i}\theta x}+\mathrm{i}\theta x\mathbf{1}%
_{(|x|<1)})\Pi ({\rm d}x)
\end{equation*}
is the so-called characteristic exponent of the process.

A property that is common to all L\'evy processes is the so-called Wiener-Hopf factorisation.  Suppose that for any $q>0$, $\mathbf{e}(q)$ is an exponentially distributed random variable with mean $q^{-1}$ that is independent of $X$. Recall that $\overline{X}_t = \sup_{s\leq t} X_s$ and let $\underline{X}_t: = \inf_{s\leq t}X_s$. The Wiener-Hopf factorisation states that the random variables $\overline{X}_{\mathbf{e}(q)}$ and $\overline{X}_{\mathbf{e}(q)}- X_{\mathbf{e}(q)}$ are independent. Thanks to the so-called principle of duality, that is to say the equality in law of the pair $\{X_{(t-s)-} - X_t: 0\leq s\leq t\}$ and $\{-X_{s}: 0\leq s\leq t\}$, it follows that $\overline{X}_{\mathbf{e}(q)}- X_{\mathbf{e}(q)}$ is equal in distribution to $-\underline{X}_{\mathbf{e}(q)}$. This leads to the following factorisation of characteristic functions
\[
\mathbb{E}({\rm e}^{{\rm i} \theta X_{\mathbf{e}(q)}}) = \mathbb{E}({\rm e}^{{\rm i}\theta\overline{X}_{\mathbf{e}(q)}})\times\mathbb{E}({\rm e}^{{\rm i}\theta \underline{X}_{\mathbf{e}(q)}}),\qquad \text{for all} \ \ \theta\in\mathbb{R}, 
\]
known as the Wiener-Hopf factorisation. Equivalently,
\begin{equation}
X_{\mathbf{e}(q)}\stackrel{d}{=}S_{q}+I_{q}, \label{prob-form}%
\end{equation}
where $S_{q}$ and $I_{q}$ are independent and equal in distribution to $\overline{X}_{\mathbf{e}(q)}$ and $\underline{X}_{\mathbf{e}(q)}$, respectively. Here we use the notation $\stackrel{d}{=}$ to mean equality in distribution.

We will work under the following assumptions:
\begin{itemize}
\item[(A1)]  $\int_{|x|\geq 1}x^2 \Pi({\rm d}x)<\infty$,

\item[(A2)] the payoff function $F:\mathbb{R}\times\mathbb{R}^{+}\to\mathbb{R}$ is a Lipschitz function with Lipschitz constant assumed to be $1$ for simplicity.

\item[(A3)] the computational time to sample from $\underline{X}_{\mathbf{e}(q)}$ and $\overline{X}_{\mathbf{e}(q)}$ is independent of the value of  $q$.
\end{itemize}
The first assumption asks for much less than what is commonly accepted in the financial literature, where typically exponential moments of the truncated L\'evy measure are required. Note that  this condition ensures  that for each $t\geq 0$, $X_t$ has a finite second (and hence also a first) moment. 

The second and third assumptions are equivalent to the corresponding conditions imposed in Dereich and  Heidenreich \cite{DH11} and Dereich \cite{D11}. A justification for the third assumption in the case that $X$ belongs to the so-called $\beta$-class of L\'evy processes is provided in Section \ref{Numerics}. As we shall discuss in more detail below, the $\beta$-class is a large family of L\'evy processes which is both widely suitable for use in financial modelling, as well as for implementation of our algorithm.

\subsection{The Wiener-Hopf Monte Carlo method}\label{sect.WHMC.method}

Due to the independent increments of L\'evy processes, the most common approaches to the Monte Carlo simulation of expectations involving the joint law of $(X_t, \overline{X}_t)$ {work with} random walk approximations to the L\'evy process. This requires one either to be able to simulate the increments of the L\'evy process exactly for a fixed time step or to be able to suitably approximate  the L\'evy process, typically by a jump diffusion process.

Recently introduced by Kuznetsov et al.~\cite{KKPvS}, the so-called Wiener-Hopf Monte Carlo method is related to the simulation of increments. However, it does not work with a fixed deterministic grid, but instead requires the underlying grid to be random with independent  and exponentially distributed spacings (i.e. the arrival times of a compound Poisson process). This allows us to simulate the paths of a very large class of L\'evy processes and thus to sample from the law of $(X_{\boldsymbol{\tau}}, \overline{X}_{\boldsymbol{\tau}})$ where $\boldsymbol{\tau}$ is a random time whose distribution can be concentrated arbitrarily close around $t$ depending on a parameter chosen in the algorithm that controls the resolution of the random grid and thus the amount of work. 

To describe the WHMC method in more detail, let us suppose that ${\mathbf{e}_{1}(1), \mathbf{e}_{2}(1), \cdots}$ are a sequence of i.i.d. exponentially distributed random variables with unit mean. The basis of the {WHMC} algorithm is the following simple observation, which follows directly from the Strong Law of Large Numbers. For all $t>0$,
\begin{equation}
\sum_{i=1}^{n} \frac{t}{n}\mathbf{e}_{i}(1) \rightarrow t \text{ as }%
n\uparrow\infty\label{SLLN}%
\end{equation}
almost surely. Note that the random variable on the left hand side of (\ref{SLLN}) can also be written as the sum of $n$ independent random variables with an exponential distribution having mean $t/n$ and therefore is equal  in law to a Gamma random variable with parameters $n$ and $n/t$, henceforth written as $\mathbf{g}(n, n/t)$. For sufficiently large $n$, Kuznetsov et al. \cite{KKPvS} argue that a suitable approximation to $\mathbb{P}(X_{t}\in\mathrm{d}x, \, \overline{X}_{t} \in\mathrm{d}y)$ is $\mathbb{P}(X_{\mathbf{g}(n,n/t)} \in\mathrm{d}x, \,\overline{X}_{\mathbf{g}(n,n/t)} \in\mathrm{d}y)$. {Indeed, it is a triviality to note that, thanks to (\ref{SLLN}) and the independence of $\mathbf{g}(n,n/t)$ from $X$, the pair $(\overline{X}_{\mathbf{g}(n,n/t)},X_{\mathbf{g}(n,n/t)} )$ converges almost surely to $(\overline{X}_t, X_t)$ as $n\uparrow\infty$}. We will  compute rates of convergence in Section \ref{num-anal}.

The following theorem is  straightforward to prove using  (\ref{prob-form}) together with the stationary and independent increments of the underlying L\'evy process. 

\begin{theorem}[Kuznetsov et al. {\cite[Thm.~1]{KKPvS}}]\label{thm_main}
Let $\{S^{j}_{n/t}: j\ge1\}$ and $\{I^{j}_{n/t}: j\ge1\}$ be i.i.d. sequences of 
random variables with common distribution equal to that of 
$\overline{X}_{\mathbf{e}(n/t)}$ and $\underline{X}_{\mathbf{e}(n/t)}$, respectively.
Then, for all $n \in\mathbb{N}$,
\[
(X_{\mathbf{g}(n,n/t)}, \overline{X}_{\mathbf{g}%
(n,n/t)}) \overset{d}{=}(V(n,n/t), J(n,n/t)),
\]
where, for any $k\in\mathbb{N}$,  and setting $V(0,n/t): = 0$ we define %and $q>0$,
\begin{align}
V(k,n/t)  &  = V(k-1,n/t) + %\sum_{j = 1}^n
\left(S^{k}_{n/t}+ I^{k}_{n/t}\right), \label{1}\\
J(k,n/t)  &  = \bigvee_{j=1}^k\left\{ V(j-1,n/t) %\sum_{k = 1}^{j-1}\left(S^{k}_{q}+ I^{k}_{q}\right)
\;+\;S^{j}_{n/t} \right\}. \label{2}%
\end{align}
\end{theorem}

It is clear from earlier remarks on the convergence of $\left({X}_{\mathbf{g}(n,n/t)},\overline{X}_{\mathbf{g}(n,n/t)} \right)$ that the pair $(V(n, n/t), J(n, n/t))$ converges in distribution to $(X_{t}, \overline{X}_{t})$. Theorem \ref{thm_main} suggests that {as soon as we are able} to simulate i.i.d. copies of the distributions of $S_{n/t}$ and $I_{n/t}$, then by the simple functional transformations given in (\ref{1}) and (\ref{2}), we may produce {an exact draw from} the  distribution of $(X_{\mathbf{g}(n, n/t)},\overline{X}_{\mathbf{g}(n, n/t)})$. Moreover, for  a suitably nice function $F$, using standard Monte Carlo methods based on the Strong Law of Large Numbers, one may estimate $\mathbb{E}(F(X_{\mathbf{g}(n,n/t) }, \overline{X}_{\mathbf{g}(n,n/t)} ) )$ by 
\begin{equation}
\widehat{F}_{\mathrm{MC}}^{n,M}:=\frac{1}{M}\sum_{i=1}^{M}
F^{n,(i)},
\label{WH-MC}%
\end{equation}
where $F^{n,(i)}$ is the $i$-th sample of 
\[
F^{n}: = F\left(V(n, n/t), J(n, n/t)\right) \,.
\]
Indeed, we have $\lim_{M\uparrow\infty} \widehat{F}_{\mathrm{MC}}^{n,M}  = \mathbb{E}\left(F(X_{\mathbf{g}(n,n/t) }, \overline{X}_{\mathbf{g}(n,n/t)} ) \right)$ almost surely, which in turn  converges   to $\mathbb{E}(F(X_{t}, \overline{X}_{t }))$ as $n\uparrow\infty$. 
  
As alluded to above, the WHMC method is numerically feasible only if samples from the distributions of $S_{n/t}$ and $I_{n/t}$ are available. Until recently, this would have proved to be a significant stumbling block on account of there being few examples for which the aforesaid distributions are known in explicit form. However, developments in Wiener-Hopf theory for L\'evy processes in the last couple of years (see for example Kuznetsov \cite{Kuz} or Kuznetsov et al. \cite{KKP2011, KKPvS}) have provided a rich enough variety of examples for which the necessary distributional sampling can be performed. This family of processes are named meromorphic L\'evy processes in Kuznetsov et al. \cite{KKP2011, KKPvS}. One large subfamily of such processes is the $\beta$-class of L\'evy processes, which also conveniently offers all the desirable properties of better known L\'evy processes that are used in mathematical finance, such as CGMY processes, VG processes or Meixner processes; see for example the discussions in  Ferreiro-Castilla and Schoutens \cite{FS11} and Schoutens and van Damme \cite{SD10}.

\subsection{The WHMC method in context}\label{WHS}

Let us now spend a little time contextualising the WHMC method for simulating the path of a L\'evy process against related possible alternative schemes, pointing out in particular where we should expect to see improved efficiencies in working with the Wiener-Hopf factorisation. A general point of reference for further reading is the book of Cont and Tankov \cite{CT}.

Generally speaking, there are three alternative approaches to simulate the path of a L\'evy process. The first one relies on the ability to simulate the increments of the L\'evy process exactly for a fixed (deterministic) time step and therefore construct a random walk as discussed at the beginning of Section \ref{sect.WHMC.method}. This method requires knowledge of the distribution of the L\'evy process at a fixed time, either through an exact analytical formula, or via numerical inversion of the characteristic function. There are also general results on the simulation of infinite divisible distributions, for which the distribution of a L\'evy process at a fixed time is an example, see e.g., Bondesson \cite{Bondesson82}. As mention earlier, approximating a L\'evy process by an embedded random walk may introduce significant errors on path functionals of $X$ such as the running maximum which is of prime interest in applications we have in mind (cf. Broadie and Glasserman~\cite{BG97}).

The second approach is to use a time-dependent infinite series expansion to approximate the value of the L\'evy process at each fixed time. A general result in this direction is found in Rosi\'nski \cite{Rosinski01} for the case that an explicit expression for the L\'evy measure  is known in closed form. The aforesaid series representation converges uniformly and almost surely in any compact set of time. This  lends itself better to sampling path functionals of the process than, perhaps the random walk approximation but it might make the numerical analysis difficult. 

Finally, the third and most common approach is to approximate $X$ via a jump-diffusion process; that is to say a L\'evy process which can be written as a linear Brownian motion plus an independent compound Poisson process. This is done by truncating the L\'evy measure, removing all small jumps below a certain threshold in magnitude and compensating for their removal by making an appropriate  adjustment to the linear and/or Gaussian component. The truncation of small jumps ensures that the remaining jumps conform to a compound Poisson structure and hence one is left with simulating the path of a linear Brownian motion, interlaced with jumps distributed according to the normalised truncated L\'evy measure arriving at an appropriate Poissonian rate. The method of truncating small jumps, originally due to Asmussen and Rosi\'nski \cite{AR01},  is an obvious approach (cf. Cont and Tankov \cite[Sect.~6.3]{CT} and the references therein) and it is the approach taken by Dereich and  Heidenreich~\cite{DH11} in their design of multilevel methods for L\'evy-driven stochastic differential equations. There it is observed that when the jump component of $X$  is of finite variation, one may reasonably replace the small jumps by a linear trend. Further, in the case of an unbounded variation jump component it is seen that a more appropriate approximation is to replace the small jumps by a Gaussian process. This truncation method is also the approach used in Dereich \cite{D11}. Whilst \cite{DH11, D11} successfully demonstrate the convergence of this truncation method, there are limitations; see, for example,  the discussion in  Asmussen and Rosi\'nski \cite{AR01}.

The main differences between these three methods and the WHMC simulation scheme are three--fold. Firstly, whilst the the WHMC scheme is restricted to L\'evy processes that allow sampling from the distributions of the variables $\overline{X}_{\mathbf{e}(q)}$ and $\underline{X}_{\mathbf{e}(q)}$, it is otherwise indifferent to the jump structure of the underlying L\'evy process. Indeed, the approximation of a L\'evy process by a compound Poisson process, which is the most popular approach, relies on the L\'evy-It\^o decomposition. The latter  decomposes the process with respect to the endogenous jump structure. In contrast, the Wiener-Hopf factorisation concerns the decomposition of the path of a L\'evy process according to the distribution of its maximum. Secondly, the WHMC method relies entirely on a randomised time grid. Note for comparison that the approach in \cite{DH11, D11} requires a deterministic grid interlaced by a random grid capturing the large jumps. Thirdly, the WHMC works on the principle of sampling the L\'evy process over the ``wrong'' time horizon (but in a way which can be made arbitrarily close to the desired time horizon) purely in order to be able to sample from the exact maximum. Note that a similar idea of randomising the time horizon in a computation involving the  expectation of a functional of the path of a linear Brownian motion that appears in the context of pricing American options is described in Carr \cite{Carr98}. We shall also see in the forthcoming numerical analysis that this idea of sampling over the ``wrong'' time horizon turns out to make the numerical analysis of the WHMC method, as well as its multilevel version, relatively straightforward. Indeed the rate of convergence of the WHMC and of the MLWH method can be expressed directly in terms of the rate of convergence of the randomised time horizon to the fixed time $t$. 

Since the Wiener-Hopf technique naturally leads to a random walk with exponential time-spacings that will on average shrink as $n \uparrow \infty$, the WHMC and the MLWH methods are in principle extendable to more general problems than the one considered here, for instance to approximate solutions of L\'evy driven SDEs (as done in \cite{DH11, D11}).

\section{Multilevel Wiener-Hopf Monte Carlo simulation}\label{MLWH_sec}

As noted already in the context of pure diffusion processes by Giles \cite{Giles07,Giles08} and others, optimal convergence rates in simulating from random processes are only achievable by using multiple levels of time grids. This has its roots in multigrid techniques for deterministic differential equations. It relies on the fact that large features or longtime trends of a process can be approximated on a coarse time grid to a sufficient accuracy that can be related clearly to the time step size. Smaller features are then added as corrections to this longtime trend %simulated on the coarse grid 
on a finer time grid. This can be done without introducing additional bias, and the benefit is that the variance of this correction on the finer grid is substantially smaller than the variance of the original process on the same grid, leading to a {\em variance reduction} and thus a smaller number of Monte Carlo samples, necessary to achieve a prescribed absolute tolerance. This variance reduction stems from the fact that our time discretisation scheme is a convergent process and, provided the convergence rates are known, the variance reduction, and thus the cost of the method, can be rigorously quantified. Applying this idea recursively on a sequence of ever finer grids, allows the simulation of the entire process with all small features in optimal computational complexity for any prescribed tolerance of the bias and the standard deviation. In more general terms, we are simulating from a random process by combining large numbers of very coarse and cheaply available paths with small numbers of fine paths that are expensive to compute. We refer to Giles \cite{Giles07,Giles08} and Cliffe et al.~\cite{CGST11} and the references therein for a detailed overview of the method, as well as to Giles and Xia \cite{GX11}, Dereich and  Heidenreich \cite{DH11} and Dereich \cite{D11} for  applications to simulating jump diffusions and L\'evy processes. Our interest here is to construct a multilevel version of the WHMC method.

\subsection{Definition of the multilevel method}\label{def_sec_MLMC}

Since the multilevel method relies on variance reduction, it is worth recalling  the formula for the mean square error of the WHMC method for subsequent reference. The mean square error between the estimator $\widehat{F}_{\mathrm{MC}}^{n,M}$ and $\E[F(X_{t},\overline{X}_{t})]$ is
\begin{eqnarray*}
e(\widehat{F}_{\mathrm{MC}}^{n,M})^{2}
&:=&\E[(\widehat{F}_{\mathrm{MC}}^{n,M}-\E[F(X_{t},\overline{X}_{t})])^{2}]  \,.
\end{eqnarray*}
Since $\E[\widehat{F}_{\mathrm{MC}}^{n,M}]=\E[F^{n}]$ and $\V[\widehat{F}_{\mathrm{MC}}^{n,M}]=M^{-1}\V[F^{n}]$, where $\mathbb{V}[\cdot]$ denotes variance,  it can be decomposed as follows:
\begin{eqnarray}
e(\widehat{F}_{\mathrm{MC}}^{n,M})^{2}&=&\E[(\widehat{F}_{\mathrm{MC}}^{n,M}-\E[\widehat{F}_{\mathrm{MC}}^{n,M}])^{2}]+(\E[\widehat{F}_{\mathrm{MC}}^{n,M}]-\E[F(X_{t},\overline{X}_{t})])^{2} \nonumber\\
&=& \V[\widehat{F}_{\mathrm{MC}}^{n,M}]+\left(\E[\widehat{F}_{\mathrm{MC}}^{n,M}]-\E[F(X_{t},\overline{X}_{t})]\right)^{2} \nonumber\\
& =&  M^{-1}\V[F^{n}]\;+\;\left(\E[F^{n}-F(X_{t},\overline{X}_{t})]\right)^{2}. \label{MSE_MC}
\end{eqnarray} 
The first term in the above decomposition is the variance of the WHMC simulation and the second one is the bias of the approximation induced by the randomised time horizon. Recall that $F^n$ converges in distribution to   $F(X_{t},\overline{X}_t)$ as $n\uparrow\infty$. We shall say that $F^n$ becomes a {\it finer} approximation of $F(X_{t},\overline{X}_t)$ the larger $n$ becomes. Conversely, the approximation $F^n$ is said to become {\it coarser} the smaller the value of $n$.

The MLWH method starts from the simple observation that we can write the expectation of the finest approximation $F^{n_{L}}$ as a telescopic sum starting from a coarser approximation $F^{n_0}$, as well as intermediate ones:
\begin{equation}\label{telescopic}
\E[F^{n_{L}}]=\E[F^{n_0}]+\sum_{\ell=1}^{L}\E[F^{n_\ell}-F^{n_{\ell-1}}]\ ,\qquad 1 \le n_0<n_1< \ldots < n_L \ .
\end{equation}
A typical choice for $n_\ell/n_0$, $\ell = 0,\ldots,L$, %(and thus also for $n/n_0$) 
are powers of $2$, or more generally of some integer $s \in \mathbb{N}$, i.e. $n_\ell := n_0  s^{\ell}$, for some $n_0 \in \mathbb{N}$. For the remainder of this paper, it suffices to simply think of $n_\ell := n_0 2^{\ell}$, for $\ell = 0,\ldots,L$. The MLWH method now consists in independently computing each of the expectations of the telescopic sum by a standard Monte-Carlo method, i.e.
\begin{equation}\label{MLMC_def}
\widehat{F}_{\mathrm{ML}}^{\mathcal{M}(n_0,L)}:=\frac{1}{M_0}\sum_{i=1}^{M_0}F^{n_0,(i)} +
\sum_{\ell=1}^{L} \frac{1}{M_{\ell}}\sum_{i=1}^{M_{\ell}}\left(F^{n_\ell,(i)}-F^{n_{\ell-1},(i)}\right)\; ,
\end{equation}
where $\mathcal{M}(n_0,L) := \{M_{\ell}\}_{\ell= 0}^{L}\,$. Analogously to (\ref{MSE_MC}), we can again expand the mean square error for the multilevel estimator to obtain
\begin{equation}\label{MSE_MLMC}
e(\widehat{F}_{\mathrm{ML}}^{\mathcal{M}(n_0,L)})^{2}:=\frac{1}{M_{0}}\V[F^{n_0}]+\sum_{\ell=1}^{L}\frac{1}{M_{\ell}}\V[F^{n_\ell}-F^{n_{\ell-1}}]\ + \ \left(\E[F^{n_{L}}-F(X_{t},\overline{X}_{t})]\right)^{2}\; .
\end{equation}
Note that the bias term remains the same, i.e. we have not introduced any additional bias. However, by a judicious choice of $\mathcal{M}(n_0,L)$ it is possible to reduce the variance of the estimator (i.e. the first two terms) for any pre--chosen computational effort, or conversely reduce the computational cost for any pre--selected tolerance of the variance of the estimator.

The variance reduction stems from the fact that $F^n$ converges to $F(X_{t},\overline{X}_{t})$ in mean square. This is equivalent to the mean square convergence of the Cauchy sequence $F^{n_\ell}-F^{n_{\ell-1}}$, and thus $\V[F^{n_\ell}-F^{n_{\ell-1}}]\rightarrow 0$. This implies that  the contribution to the total variance from the finer levels (where the cost per sample is high) is decreasing with
$n_\ell \to \infty$, allowing us to use only very few of these expensive samples, while on the coarser levels the cost of samples is smaller anyway. At this point in time, we have not yet established rigorously that $F^n$ converges to $F(X_{t},\overline{X}_{t})$, not to mention at what rate. This will be dealt with in Section \ref{num-anal}. However, let us assume mean square convergence for now.

\subsection{Thinning a Poisson random grid}\label{refinement}

Before continuing, it is first necessary to elaborate on another important point. One of the defining features of the multilevel method is that both approximations of the payoff function in simulating from the random variable $(F^{n_\ell}-F^{n_{\ell-1}})$ should come from the same draw. This means that once the increments of the L\'evy process have been sampled to obtain a draw for $F^{n_\ell}$, a deterministic transformation of the sample points and of the increments should take place to obtain a sample for $F^{n_{\ell-1}}$ (or vice versa). While this procedure is clear in the original setting of Giles \cite{Giles08}, where a diffusion processes is approximated by a random walk on a deterministic grid, it is not entirely trivial how to proceed in the case of a general L\'evy process and with respect to random grids. 

Our approach to construct a sample for a coarser approximation from the draw of a finer approximation relies on the technique of Poisson thinning. A similar approach was considered in Giles and Xia \cite{GX11} in the case of jump diffusion processes and can be traced back further (also for general L\'evy processes) to, e.g., Glasserman and Merener \cite{GM03}. 

The first step to sample from $(F^{n_\ell}-F^{n_{\ell-1}})$ on level $\ell$ is to construct a Poisson random grid with rate $n_\ell/t$, \ie the inter-arrival times are independently distributed as an exponential random variable with mean $t/n_\ell$. Fix  $\ell \in\N$ and  $t>0$ and let $N^{\ell}  :=\{N^{\ell}_s: s\geq 0\}$ denote a Poisson process with arrival rate $n_\ell/t$; the arrival times in $N^{\ell}$ are denoted by $\{T^{N,\ell}_k : k\geq 0\}$, with $T^{N,\ell}_0 = 0$. Recall that, without loss of generality, we restricted to the case where $n_\ell/n_{\ell-1} = 2$ and notice that under this notation we have $\mathbf{g}(n_\ell,n_\ell/t)=T^{N,\ell}_{n_\ell}$. It is a well established fact (cf. Kingman \cite{king}) that if we censor arrivals in the process $N^\ell$ by tossing independent fair coins at each arrival and ignoring the arrival if the coin lands, say, on heads, then the resulting point process in time has the same law as $N^{\ell-1}$. In particular the inter-arrival times in this coarser grid are independent and exponentially distributed with mean $t/n_{\ell-1}$. The general case $n_\ell/n_{\ell-1} = s$, for some $s>1$, could clearly be treated in an analogous manner by tossing a biased coin with probability of acceptance $1/s$.

To make the construction more precise, let us return to the case $s=2$ and consider the sequences  $\{S^{j}_{n_{\ell}/t}: j\geq1\}$ and $\{I^{j}_{n_{\ell}/t}: j\geq1\}$ of i.i.d. random variables with common distributions $\overline{X}_{\mathbf{e}(n_{\ell}/t)}$ and  $\underline{X}_{\mathbf{e}(n_{\ell}/t)}$ respectively, where the exponential periods are taken from the Poisson process $N^\ell$. Suppose we set $\kappa_0 =0$,  and for $i\geq 1$, we write $\kappa_i$ for the index of the $i$-th accepted arrival in the process $N^\ell$ after censoring. Note that for $i\geq 1$, $\kappa_i - \kappa_{i-1}$ are a sequence of i.i.d. geometrically distributed random variables with parameter one half. Note also that, for $i\geq 1$, the time that elapses in the process $N^\ell$ between the $(i-1)$-th (after censoring) arrival and the $i$-th (after censoring) arrival is equal in distribution to precisely
\begin{equation}
\sum_{j = \kappa_{i-1}+1}^{\kappa_i} \mathbf{e}_j(n_\ell/t).
\label{geom-exp}
\end{equation}
It is a straightforward computation to show that the above random variable is exponentially distributed with mean $t/n_{\ell-1}$, e.g. by computing its moment generating function.  

Thanks to the Poisson thinning, we can now construct the independent sequences $\{S^{i}_{n_{\ell-1}/t}: i\geq1\}$ and $\{I^{i}_{n_{\ell-1}/t}: i\geq1\}$ of random variables, which are needed for sampling from $F^{n_{\ell-1}}$, via a deterministic transformation of the sequences  $\{S^{j}_{n_{\ell}/t}: j\geq1\}$ and $\{I^{j}_{n_{\ell}/t}: j\geq1\}$ as follows:
\begin{eqnarray}
S^{i}_{n_{\ell-1}/t}&=&\bigvee_{k=1}^{\kappa_i - \kappa_{i-1}} \left\{\sum_{j=1}^{k-1}\left(S_{n_{\ell}/t}^{\kappa_{i-1}+j}+I_{n_{\ell}/t}^{\kappa_{i-1}+j}\right) \; + \; S_{n_{\ell}/t}^{\kappa_{i-1}+k} %:\ k=1, \cdots, (\kappa_i - \kappa_{i-1})
\right\}
\label{refinement_1}\\
I^{i}_{n_{\ell-1}/t}&=&\sum_{j=\kappa_{i-1}+1}^{\kappa_i}\left(S_{n_{\ell}/t}^{j}+I_{n_{\ell}/t}^{j}\right)\;-\;S^{i}_{n_{\ell-1}/t}\label{refinement_2}
\end{eqnarray} 
Although it is now clear how to construct the infinite sequence of pairs $\{(S^{i}_{n_{\ell-1}/t}, I^{i}_{n_{\ell-1}/t}) : i\geq 1\}$ from $\{(S^{j}_{n_{\ell}/t}, I^{j}_{n_{\ell}/t}): j\geq 1\}$, some explanation is still needed when dealing with finite sequences, such as those required by the MLWH method. When the Poisson random grid is stopped at the time horizon $\mathbf{g}(n_{\ell}, n_{\ell}/t)$, it will in general be necessary to first extend the number of exponential periods from the process $N^\ell$ beyond $\mathbf{g}(n_\ell, n_\ell/t)$ before thinning to get the corresponding $\ell-1$ level grid. This is because we cannot ensure that after thinning, the total non-censored points in $N^{\ell}$ up to $ \mathbf{g}(n_\ell, n_\ell/t)$ sum up to $n_{\ell-1}$. Alternatively, one can produce the remaining terms of the coarser level by sampling directly from $(S_{n_{\ell-1}/t}, I_{n_{\ell-1}/t})$ due to the thinning theorem.

Let us point out here another issue directly related to the fact that we are constructing our algorithm based on a completely random grid. Denote by $\check{F}^{n_{\ell-1}}$ a sample produced from $F^{n_{\ell}}$ by the thinning methodology described above. Since the telescopic sum in (\ref{telescopic}) has to cancel out, \ie we do not want to introduce extra bias, it is important that the expectation of $F^{n_{\ell-1}}$ is the same as the expectation of $\check{F}^{n_{\ell-1}}$. This is called the consistency of the multilevel algorithm. Thanks to the thinning theorem this is straight forward to check in the above construction, however it will play a role in Section \ref{randomisation_time}. Indeed, the above description asserts that thinning a Poisson process $N^{\ell}$ will lead in distribution to a Poisson process of rate $n_{\ell-1}/t$ and is thus equal in distribution to $N^{\ell-1}$. Therefore the random time $\mathbf{g}(n_{\ell-1}, n_{\ell-1}/t)$ has the same distribution as the resulting thinned version from $\mathbf{g}(n_{\ell}, n_{\ell}/t)$, say $\check{\mathbf{g}}(n_{\ell-1}, n_{\ell-1}/t)$. Finally, due to the fact that $\mathbf{g}(n_{\ell-1}, n_{\ell-1}/t)$ and $\check{\mathbf{g}}(n_{\ell-1}, n_{\ell-1}/t)$ are independent of the underlying L\'evy process the consistency follows as $\check{F}^{n_{\ell-1}}$ and ${F}^{n_{\ell-1}}$ have the same law.

\section{Analysis of Wiener-Hopf Monte Carlo simulation} \label{num-anal}

Henceforth we shall use the following notation. We will write $a\lesssim b$ for two positive quantities $a$ and $b$, if $a/b$ is uniformly bounded independent of any parameters, such as $t$, $L$, $\{n_\ell\}_{\ell=0}^L$, or $\{M_\ell\}_{\ell=0}^L$, as well as of any parameters describing the underlying L\'evy process. We will write $a\eqsim b$, if $a\lesssim b$ and $b\lesssim a$. 

\subsection{Abstract convergence theorems}

There are two ways to quantify and compare the complexity of algorithms. We can either specify a target accuracy $\varepsilon$ and then estimate the cost to achieve an error below this target accuracy with a particular method, or conversely, we can analyse the convergence of an algorithm as a function of the cost. The results are easily interchangeable. The first approach is the one chosen in Giles \cite{Giles08} and Cliffe et al.~\cite{CGST11}. However, Dereich and  Heidenreich \cite{DH11} and Dereich \cite{D11} chose the second approach and since we will mainly compare our results with theirs, we will also choose this approach here. To fix notation, computational cost is measured in floating point operations here (but it could equally be CPU-time). Thus, for the following let us assume that the expected value of the computational cost $\mathcal{C}(\widehat{F})$ to compute a particular estimator $\widehat{F}$ of $\E[F(X_{t},\overline{X}_{t})]$ is bounded by~$\nu$, i.e. 
$$\mathbb{E} \big[\mathcal{C}(\widehat{F})\big] \lesssim \nu\,.$$
Then we want to bound the root mean square error $e(\widehat{F})$ in terms of $\nu$.

Let us start by analysing the single--level WHMC method. The following result can be easily deduced from (\ref{MSE_MC}) with the modification that the cost per sample is itself a random variable. It will become clear in the proof of Theorems \ref{g-error} and \ref{T-error} why we need this generalization.

\begin{theorem}
\label{thm_WHMC}
Let $t>0$ and let us assume that $F^{n}$ converges in mean square to $F(X_{t},\overline{X}_{t})$. Suppose further that there exist positive constants $\alpha, \gamma > 0$ such that
\begin{enumerate}
\item[(i)] $|\E[F^{n}-F(X_{t},\overline{X}_{t})]|\lesssim n^{-\alpha}$,
\item[(ii)] $\E[\mathcal{C}_{n}]\lesssim n^{\gamma}$,
\end{enumerate}
where $\mathcal{C}_{n}$ represents the cost of computing a single sample of $F^{n}$ on a Poisson grid with rate~$n/t$. Then, for every $\nu \in \mathbb{N}$, there exist $n, M \in \mathbb{N}$ such that 
\[
\E\left[\mathcal{C}\left(\widehat{F}_{\mathrm{MC}}^{n,M}\right)\right] \lesssim \nu \qquad \text{and} \qquad e\left(\widehat{F}_{\mathrm{MC}}^{n,M}\right) \lesssim \nu^{-\frac{1}{2+\gamma/\alpha}}\,.
\] 
\end{theorem}
\proof
Assuming for the moment that $\mathbb{V}[F^{n}]$ is bounded independently of $n$, then balancing the two terms on the right hand side of (\ref{MSE_MC}) and using assumption (i), we see that we should choose $M \eqsim n^{2\alpha}$. Now, since
\[
\mathcal{C}\left(\widehat{F}_{\mathrm{MC}}^{n,M}\right) = M \, \mathcal{C}_{n},
\]
we can deduce from assumption (ii) that the expected value of the total cost of the estimator will be bounded by $\nu$, if we choose $n \eqsim \nu^{\frac{1}{2\alpha+\gamma}}$ and $M \eqsim \nu^{\frac{2 \alpha}{2\alpha+\gamma}}$, which leads to the required bound for $e\big(\widehat{F}_{\mathrm{MC}}^{n,M}\big)$. 

It remains to bound $\mathbb{V}[F^{n}]$. First note that trivially $\mathbb{V}[F^{n}] \le \mathbb{E}[(F^{n})^2]$ and so 
\begin{equation}\label{var_cte}
\frac{1}{2} \,\mathbb{V}[F^{n}] \le \mathbb{E}[(F^{n} - F(X_{t},\overline{X}_{t}))^2]
\;+\; \mathbb{E}[(F(X_{t},\overline{X}_{t})^2]
\end{equation}
The first term can be bounded independently of $n$ using the Lipschitz continuity of $F$ in Assumption (A2) and our assumption that $F^{n}$ converges in mean square to $F(X_{t},\overline{X}_{t})$. The second term is bounded due to the Lipschitz continuity of $F$ and Assumption (A1). 
\qed

Thanks to the general considerations in Giles \cite{Giles08} and Cliffe et al. \cite[Theorem 1]{CGST11}, the analysis of the multilevel version follows in the same way. The proof of the following result is similar to that in \cite[Appendix A]{CGST11}, with the same modifications as above to deal with random computational cost.
\begin{theorem}
\label{Rob_thm}
Let $t>0$ and $n_\ell = n_0 s^\ell$, for some $\ell,n_0 \in \mathbb{N}$ and $s>1$, and 
suppose that there are positive constants $\alpha, \beta,\gamma > 0$ with
$\alpha\geq\frac{1}{2}(\beta\wedge\gamma)$ such that
\begin{enumerate}
\item[(i)] $|\E[F^{n_\ell}-F(X_{t},\overline{X}_{t})]|\lesssim n_\ell^{-\alpha}$
\item[(ii)] $\V(F^{n_\ell}-F^{n_{\ell-1}})\lesssim n_\ell^{-\beta}$
\item[(iii)] $\E[\mathcal{C}_{n_\ell}]\lesssim n_\ell^{\gamma}$,
\end{enumerate}
where $\mathcal{C}_{n}$ represents the cost of computing a single sample of $F^{n}$ on a Poisson grid with rate~$n/t$. Then, for every $\nu\in \mathbb{N}$, there exists a value $L$ and a sequence $\mathcal{M}(n_0,L) = \{M_{\ell}\}_{\ell=0}^L$ such that%\vspace{-2ex}
\begin{equation*}
\E\left[\mathcal{C}\left(\widehat{F}_{\mathrm{ML}}^{\mathcal{M}(n_0,L)}\right)\right] \lesssim \nu \qquad \text{and} \qquad 
e\left(\widehat{F}_{\mathrm{ML}}^{\mathcal{M}(n_0,L)}\right) 
%:= \sqrt{\E[(\widehat{F}_{\mathrm{ML}}^{\mathcal{M}(n_0,L)}-\E[F(X_{t},\overline{X}_{t})])^{2}]}
\lesssim \left\{
\begin{array}{ll}
\nu^{-\frac{1}{2}}\,,&\text{if }\ \beta>\gamma\,,\\
\nu^{-\frac{1}{2}}\log^2\nu\,, &\text{if }\ \beta=\gamma\,,\\
\nu^{-\frac{1}{2+(\gamma-\beta)/\alpha}}\,, &\text{if } \ \beta<\gamma\,.
\end{array}
\right.
\end{equation*}
\end{theorem}

In order to apply Theorems \ref{thm_WHMC} and \ref{Rob_thm} now, we need to find suitable $\alpha,\; \beta$ and $\gamma$ such that the necessary assumptions are satisfied. We will verify in the next subsection that the values of $\alpha$ and $\beta$ in the hypothesis of Theorems \ref{thm_WHMC} and \ref{Rob_thm} are directly related to the convergence of the time horizon $\mathbf{g}(n,n/t)$ to $t$.
As for the value of $\gamma$, Assumption (A3) in Section~\ref{WHMC_MLMC} implies that for both the single and the multilevel algorithm we always have $\gamma = 1$. This is clear for the single level case, since the cost is not random, but a proof is needed for the multilevel version as we might have to enlarge the finer grid in the thinning procedure as described in Section \ref{refinement}. 
\begin{lemma}\label{cost_gamma_1}
Under the notation of Theorem \ref{Rob_thm}, assume that $\mathcal{C}(F^{n_\ell})\lesssim n_\ell$. Then the expected cost of producing a sample of $F^{n_\ell}-F^{n_{\ell-1}}$ via the thinning methodology of Section \ref{refinement} is 
\begin{equation*}
\E[\mathcal{C}(F^{n_\ell}-F^{n_{\ell-1}})]\lesssim n_\ell\ .
\end{equation*}
\end{lemma}
\proof
Let us first focus on the case that $n_\ell/n_{\ell-1}=2$ and the thinning methodology corresponds to tossing a fair coin. In this case, the probability of choosing $i$ arrival times from the finer level when thinning is $\binom{n_\ell}{i}{2^{-n_\ell}}$ and hence one easily concludes from the assumptions that
\begin{eqnarray*}
\E[\mathcal{C}(F^{n_\ell}-F^{n_{\ell-1}})] 
&\lesssim& n_\ell+ \sum_{i=0}^{n_\ell/2}\binom{n_\ell}{i}\frac{1}{2^{n_\ell}}\mathcal{C}(F^{\frac{n_\ell}{2}-i})\\
&\lesssim& n_\ell+ \sum_{i=0}^{n_\ell/2}\binom{n_\ell}{i}\frac{1}{2^{n_\ell}}\left(\frac{n_\ell}{2}-i\right)\\
&=&n_\ell+\frac{\Gamma(n_\ell)}{(\Gamma(n_\ell/2))^{2}2^{n_\ell}}\\
&\simeq&n_\ell\ ,
\end{eqnarray*} 
where $\mathcal{C}(F^{\frac{n_\ell}{2}-i})$ corresponds to the cost of producing the remaining terms on the coarser level after choosing $i$ terms on the finer level. 

In the general case where $n_\ell/n_{\ell-1}=s$ the probability \mbox{$\binom{n_\ell}{i}{2^{-n_\ell}}$} should be replace by \mbox{$\binom{n_\ell}{i}s^{i}(1-s)^{n_\ell-i}$} and a similar result holds, although the computations are more tedious. \qed

\subsection{Explicit convergence rates}
\label{g_time} 
 
Let us write for convenience $\mathbf{X}_t=(X_t,\overline{X}_t)$, for $t\geq 0$, and let $\boldsymbol{\tau}$ and $\boldsymbol{\tau}'$ be any non-negative random variables that are independent of $X$ (but potentially correlated). Then, thanks to the Lipschitz assumption (A2) on $F$, it is straightforward to deduce that % for any $t>0$,
\begin{eqnarray}
\label{eq_alpha}
|\E[F(\mathbf{X}_{\boldsymbol{\tau}})-F(\mathbf{X}_{t})]|
&\leq&\E[|\mathbf{X}_{\boldsymbol{\tau}}-\mathbf{X}_{t}|]\ \leq\ \E[|X_{\boldsymbol{\tau}}-X_{t}|] \;+\; \E[|\overline{X}_{\boldsymbol{\tau}}-\overline{X}_{t}|]
\end{eqnarray}
and 
\begin{eqnarray}
\label{eq_beta}
\V(F(\mathbf{X}_{\boldsymbol{\tau}})-F(\mathbf{X}_{\boldsymbol{\tau}'}))
&\leq& \E[(\mathbf{X}_{\boldsymbol{\tau}}-\mathbf{X}_{\boldsymbol{\tau}'})^{2}] %\\
\ \leq\  2 \Big(\E[(X_{\boldsymbol{\tau}}-X_{\boldsymbol{\tau}'})^{2}] \;+\; \E[(\overline{X}_{\boldsymbol{\tau}}-\overline{X}_{\boldsymbol{\tau}'})^{2}]\Big). %\notag
\end{eqnarray}
To estimate the quantities on the right hand sides of (\ref{eq_alpha}) and (\ref{eq_beta}) we will make use of the following lemma.
\begin{lemma}\label{mean_square_error_general}
Suppose that  $X$ is a L\'evy process which satisfies Assumption (A1). Let us denote by $\boldsymbol{\tau}$ any non-negative random variable, independent of $X$. Then, for any fixed $t>0$,
\begin{enumerate}
\item[(i)] $\E[(X_{\boldsymbol{\tau}}-X_{t})^{2}]\;=\;\V[X_{1}]\,\E[|\boldsymbol{\tau}-t|]\,+\,(\E[X_{1}])^{2}\,\E[(\boldsymbol{\tau}-t)^2]$,
\item[(ii)] $\E[(\overline{X}_{\boldsymbol{\tau}}-\overline{X}_{t})^{2}]\;\leq \; 16 \V[X_{1}]\,\E[|\boldsymbol{\tau}-t|]\,+\,2(\E[X_{1}] \vee 0)^{2})\,\E[(\boldsymbol{\tau}-t)^2]$.
\end{enumerate}
\end{lemma}
\proof
(i) \;Due to the fact that increments of a L\'evy process are i.i.d. we can write
\begin{equation*}
X_{\boldsymbol{\tau}}-X_{t}\stackrel{d}{=}
\left\{
\begin{array}{cr}
X_{\boldsymbol{\tau}-t}\,,&\text{if } \boldsymbol{\tau}\geq t,\\
-X_{t-\boldsymbol{\tau}}\,,&\text{if } \boldsymbol{\tau}<t.
\end{array}
\right.
\end{equation*}
Since $X$ and $\boldsymbol{\tau}$ are assumed to be independent and since $\mathbb{E}[X_s^2] =\mathbb{V}[X_1]\, s +\mathbb{E}[X_1]\,s^2 $, for all $s\geq 0$, this implies
\begin{eqnarray*}
\E[(X_{\boldsymbol{\tau}}-X_{t})^{2}]
&=&\V[X_{1}]\E[|\boldsymbol{\tau}-t|]+(\E[X_{1}])^{2}\E[(\boldsymbol{\tau}-t)^2]\;.
\end{eqnarray*}

\noindent
(ii) \;To obtain the result for the mean square error of the supremum process we first note the following equality in distribution
\begin{equation*}
\overline{X}_{t}\stackrel{d}{\;=\;}\overline{X}_{s} \,\vee\, (X_{s}+\overline{X}'_{t-s}),
\end{equation*}
where $s<t$ and $\overline{X}'_{t-s}$ is independent of $\{X_{s'}: s'\leq s\}$ and  identically distributed to  $\overline{X}_{t-s}$. Taking account of the duality property for L\'evy processes,  which tells us that $\underline{X}_{s}$ is equal in distribution to $X_s  -\overline{X}_s$, we have that
\begin{equation*}
\overline{X}_{t}-\overline{X}_{s}
\stackrel{d}{\;=\;}0\,\vee\, (\underline{X}''_{s}+\overline{X}'_{t-s}) \;\leq\; \overline{X}'_{t-s},
\label{hat-diff}
\end{equation*}
where $\underline{X}''_{s}$ is independent of $\overline{X}'_{t-s}$ and equal in distribution to   $\underline{X}_{s}$. 
It is now easy to check that 
\begin{equation*}
\E[(\overline{X}_{\boldsymbol{\tau}}-\overline{X}_{t})^{2}]\leq
\E[\overline{X}_{|\boldsymbol{\tau}-t|}^{2}].
\end{equation*}
Let us decompose $X$ into its martingale part $X^{*}:=\{X^{*}_{t}:t\geq0\}$ and a drift, \ie $X_{t}=X^{*}_{t}+t\E[X_1]$, for $t\geq 0$. Due to Assumption (A1) this is always possible. Then we have
\begin{eqnarray}
\E[\overline{X}_{t}^{2}]&=&
\E[\sup_{s\leq t}(X^{*}_{s}+\E[X_1]s)^{2}]\notag\\
&\leq&
\E[(\overline{X}^{*}_{t}+t(\E[X_1]\vee 0))^{2}]\notag\\
&\leq&
2\bigg(
\E[(\overline{X}^{*}_{t})^{2}]+(\E[X_1]\vee 0)^{2}t^2
\bigg),
\label{return1}
\end{eqnarray}
where $\overline{X}^*_t = \sup_{s\leq t }X^*_t$.
Appealing to Doob's  martingale inequality (see e.g.~Sato \cite[p.~167]{Sato99}) we end up with the inequality
\begin{equation}
\E[\overline{X}_{t}^{2}]
\leq
2\bigg(
8\E[(X^{*}_{t})^{2}]+(\E[X_1]\vee 0)^{2}\,t^{2}
\bigg),
\label{return2}
\end{equation}
from which we may now write (again remembering that $\boldsymbol{\tau}$ is independent of $X$)
\begin{equation*}
\E[\overline{X}_{|\boldsymbol{\tau}-t|}^{2}]
\leq
16\E[(X^{*}_{|\boldsymbol{\tau}-t|})^{2}]+2(\E[X_1]\vee 0)^{2}\,\E[(\boldsymbol{\tau}-t)^{2}].
\end{equation*}
Since $E[(X^{*}_{|\boldsymbol{\tau}-t|})^{2}]=\E[(X^{*}_{\boldsymbol{\tau}}-X^{*}_{t})^{2}]$, we may use the conclusion of  part (i) and finally write
\begin{eqnarray*}
\E[\overline{X}_{|\boldsymbol{\tau}-t|}^{2}]
&\leq&
{16}\E[(X^{*}_{\boldsymbol{\tau}}-X^{*}_{t})^{2}]+{2}(\E[X_1]\vee 0)^{2}\,\E[(\boldsymbol{\tau}-t)^{2}]\notag\\
&=&
{16}\V[X_{1}^{*}]\,\E[|\boldsymbol{\tau}-t|]+2(\E[X_1]\vee 0)^{2})\,\E[(\boldsymbol{\tau}-t)^2],
\label{spirit}
\end{eqnarray*}
which concludes the proof since $\V[X_{1}^{*}]=\V[X_{1}]$.
\qed

Using this lemma we can now verify the remaining hypotheses in Theorems \ref{thm_WHMC} and \ref{Rob_thm} and in particular establish the values for the parameters $\alpha$ and $\beta$ in the convergence bounds for the WHMC and MLWH methods.

\begin{prop}\label{mean_square_error_gamma}
Let $t>0$ and $X$ satisfies Assumption (A1). Then, for any $n \in \mathbb{N}$, we have 
\begin{enumerate}
\item[(i)] $\E[(X_{\mathbf{g}(n,n/t)}-X_{t})^{2}]\eqsim n^{-1/2}$ \ and \ $\E[|X_{\mathbf{g}(n,n/t)}-X_{t}|]\lesssim n^{-1/4},$
\item[(ii)] $\E[(\overline{X}_{\mathbf{g}(n,n/t)}-\overline{X}_{t})^{2}]\lesssim n^{-1/2}$ \ and \ $\E[|\overline{X}_{\mathbf{g}(n,n/t)}-\overline{X}_{t}|]\lesssim n^{-1/4}.$
\end{enumerate}
If, in addition, $X$ has paths of bounded variation then we have the sharper bound
\begin{equation*}
\E[|X_{\mathbf{g}(n,n/t)}-X_{t}|]\lesssim n^{-1/2}
\quad\text{and}\quad
\E[|\overline{X}_{\mathbf{g}(n,n/t)}-\overline{X}_{t}|]\lesssim n^{-1/2}.
\end{equation*}
\end{prop}
\proof 
It suffices to prove the results for the second moments in (i) and (ii). The results for the first moments then follow from Jensen's inequality. With a view to applying Lemma \ref{mean_square_error_general}, we will  show that for all $n\in\mathbb{N}$
\begin{eqnarray}
\E[(\mathbf{g}(n,n/t)-t)^{2}]&=&\frac{t^{2}}{n}\label{gamma_1}\\
\E[|\mathbf{g}(n,n/t)-t|]&=&2t{\mathrm e}^{-n}\frac{n^{n}}{n!} \ \eqsim \ n^{-\frac{1}{2}}\label{gamma_2}\ ,
\end{eqnarray}
where the final equivalence in (\ref{gamma_2}) follows from Stirling's formula. 

First note that, since $\E[\mathbf{g}(n,n/t)]=t$, equation (\ref{gamma_1}) is just the formula for the variance of the gamma distribution $\mathbf{g}(n,n/t)$.
For equation (\ref{gamma_2}), it is clear that
\begin{eqnarray}
\E[|\mathbf{g}(n,n/t)-t|]
&=&\int_{0}^{t}(t-s)g(n,n/t,s){\rm d} s+\int_{t}^{\infty}(s-t)g(n,n/t,s){\rm d}s\notag\\
&=&2\int_{0}^{t}(t-s)g(n,n/t,s){\rm d}s,\label{gamma_sym}
\end{eqnarray}
where $g(k,\theta,s)=\Gamma(k)^{-1}\theta^{k}s^{k-1}{\mathrm e}^{-\theta s}$.
Recall the definition of the incomplete gamma function, $\gamma(k,u):=\int_{0}^{u}x^{k-1}{\mathrm e}^{-x}{\rm d}x$, so that
\begin{equation*}
\int_{0}^{t}g(k,\theta,s){\rm d}s=\frac{\gamma(k,t\theta)}{\Gamma(k)}
\quad\text{ and }\quad
\int_{0}^{t}sg(k,\theta,s){\rm d}s=\frac{\gamma(k+1,t\theta)}{\Gamma(k)\theta}\ .
\end{equation*}
We can now rewrite the expression (\ref{gamma_sym}) as
\begin{eqnarray}
\int_{0}^{t}(t-s)g(n,n/t,s){\rm d}s
&=&t\int_{0}^{t}g(n,n/t,s){\rm d}s-\int_{0}^{t}sg(n,n/t,s){\rm d}s \notag\\
&=&t\left(
\frac{\gamma(n,n)}{\Gamma(n)}-\frac{\gamma(n+1,n)}{\Gamma(n+1)}
\right)\label{gamma_abra}
\end{eqnarray}
From the representation of $\gamma(k,u)$ in, for example,  Abramowitz and Stegun \cite[Section 6.5.13]{AS70}), we can write for $k = 1,2,\cdots$,
\begin{equation*}
\frac{\gamma(k,u)}{\Gamma(k)} =
\frac{1}{\Gamma(k)}\int_{0}^{u}x^{k-1}{\mathrm e}^{-x}{\rm d}x
=1-\left(1+u+\frac{u^{2}}{2!}+\frac{u^{3}}{3!}+\cdots+\frac{u^{k-1}}{(k-1)!}\right){\mathrm e}^{-u},
\end{equation*}
which, when substituted into (\ref{gamma_abra}), 
together with (\ref{gamma_sym}) gives the equality (\ref{gamma_2}).

If $X$ has paths of bounded variation then we can return to the proof of 
Lemma \ref{mean_square_error_general}, in particular the estimates in (\ref{return1}) and (\ref{return2}), and use similar computations together with Doob's martingale inequality for absolute moments (cf.  Sato \cite[p. 167]{Sato99}) to 
deduce that 
\begin{eqnarray}
\E[|\overline{X}_{\mathbf{g}(n,n/t)}-\overline{X}_{t}|]%&\leq&\E[\overline{X}_{|\mathbf{g}(n,n/t)-t|}]\notag\\
&\leq&\E[\overline{X}^{*}_{|\mathbf{g}(n,n/t)-t|}]+(\E[X_1]\vee 0)\E[|\mathbf{g}(n,n/t)-t|]\notag\\
&\leq&
8\E[|X^{*}_{|\mathbf{g}(n,n/t)-t|}|]+(\E[X_1]\vee 0)\E[|\mathbf{g}(n,n/t)-t|].
\label{inturn}
\end{eqnarray}
Note also that $\mathbb{E}[|{X}_{\mathbf{g}(n,n/t)}-{X}_{t}|] = \mathbb{E}[|X_{|\mathbf{g}(n,n/t)-t|}|]$ is essentially bounded by a linear combination of $\E[|X^{*}_{|\mathbf{g}(n,n/t)-t|}|]$ and $\E[|\mathbf{g}(n,n/t)-t|]$.
Taking account of (\ref{gamma_2}), we therefore see that the proof is complete 
as soon as we show that  $\E[|X^{*}_{|\mathbf{g}(n,n/t)-t|}|]\lesssim n^{-1/2}$. 

To this end, note that every L\'evy process of bounded variation can be written as the difference of two independent subordinators. Accordingly we shall write, for any $t\geq 0$, $X^{*}_{t}=X'_{t}-X''_{t}$\,. Assumption (A1) ensures that both $X'$ and $X''$ have finite first and second %, and hence first, 
moments at all fixed times. Therefore, using a similar derivation as above, as well as (\ref{gamma_2}), we have 
 \begin{eqnarray*}
 \E[|X^{*}_{|\mathbf{g}(n,n/t)-t|}|]&\leq&\E[X'_{|\mathbf{g}(n,n/t)-t|}]+\E[X''_{|\mathbf{g}(n,n/t)-t|}]\\
 &=&\E[|\mathbf{g}(n,n/t)-t|]\mathbb{E}[|X'_1| +|X''_1|] \ \eqsim \ n^{-1/2}\,.
\end{eqnarray*}
\qed

We are now ready to state the main result of the paper which gives the convergence rates of the WHMC and the MLWH method.  

\begin{teo}\label{g-error} 
Suppose that assumptions (A1)--(A3) are satisfied. Then,  the hypotheses in Theorems \ref{thm_WHMC} and \ref{Rob_thm} hold with $\alpha = 1/4$, $\beta = 1/2$ and $\gamma=1$. In the case that $X$ has paths of bounded variation we have $\alpha = 1/2$. Thus, for any $\nu \in \mathbb{N}$ and under the constraint that the expected value of the total computational cost is $\mathcal{O}(\nu)$ operations, the root mean square errors for the single and multilevel methods are
\begin{equation}\label{rate_g_WHMC}
e(\widehat{F}_{\mathrm{MC}}^{n,M})\; \lesssim \; \left\{
\begin{array}{ll}
\nu^{-\frac{1}{6}}, \ & \text{if} \ \ \text{$X$ is of unbounded variation,} \\[1ex]
\nu^{-\frac{1}{4}}, \ & \text{if} \ \ \text{$X$ is of bounded variation.}
\end{array} \right. 
\end{equation}
and
\begin{equation}\label{rate_g_MLWH}
e(\widehat{F}_{\mathrm{ML}}^{\mathcal{M}(n_0,L)})\; \lesssim \; \left\{
\begin{array}{ll}
\nu^{-\frac{1}{4}}, \ & \text{if} \ \ \text{$X$ is of unbounded variation,} \\[1ex]
\nu^{-\frac{1}{3}}, \ & \text{if} \ \ \text{$X$ is of bounded variation.}
\end{array} \right. 
\end{equation}
\end{teo}
\proof This follows immediately from Theorems \ref{thm_WHMC} and \ref{Rob_thm}, Proposition \ref{mean_square_error_gamma}, equations (\ref{eq_alpha}) and (\ref{eq_beta}) and the fact that $\gamma = 1$ (cf.~Assumption (A3) and Lemma \ref{cost_gamma_1}). 
\qed

\begin{remark}\label{weak_vs_strong_error}\rm
Note that Assumption (A2) implicitly transforms the weak error estimate in Theorems \ref{thm_WHMC} and \ref{Rob_thm} into a strong error estimate through inequality (\ref{eq_alpha}). This means that by using Proposition \ref{mean_square_error_gamma} we only get a lower bound estimate of the asymptotic parameter $\alpha$. By noting that the time horizon $\mathbf{g}(n,n/t)$ is an unbiased estimator of $t$ one easily sees that $|\E[X_{\mathbf{g}(n,n/t)}-X_{t}]|=0$. Although this does not hold for the supremum process in general, for particular cases and under regularity constrains on the function $F$, estimates for the weak convergence\footnote{Weak convergence refers here to the behaviour of $|\E[F(\mathbf{X}_{\mathbf{g}(n,n/t)})]-\E[F(\mathbf{X}_{t})]|$ as $n\uparrow\infty$.} might be computable, see for instance the approach taken in Jacod et al.~\cite{JKMP05} in the case of functions depending on the terminal value of the solution of a SDE. We have chosen to follow the approach taken in Dereich \cite{D11} and Dereich and  Heidenreich \cite{DH11} which allows for more general functions $F$ with the penalization that we will have only strong estimates. We will compute numerically the weak error in the numerical implementation of the algorithm in Section \ref{experiments} and show that there is a substantial improvement of the estimates of parameter $\alpha$.
\end{remark}

We finish this section by comparing the convergence results of Theorem \ref{g-error} with those of Dereich and Heidenreich \cite{DH11}, Dereich \cite{D11} and  Jacod et al.~\cite{JKMP05}. Since all these alternative methods show a significant dependence on the jump structure of the underlying L\'evy process, let us first recall the definition of the Blumenthal--Getoor index,
\begin{equation*}
\rho:=\inf\left\{\beta>0 : \int_{(-1,1)}|x|^{\beta}\nu({\rm d}x)<\infty\right\}\in[0,2]\ .
\end{equation*}
The higher this index, the more important small jumps are. L\'evy processes with finite activity will always have Blumenthal--Getoor index $\rho=0$.

\begin{figure}[htb]
\center
\includegraphics[height=7cm,width=9cm]{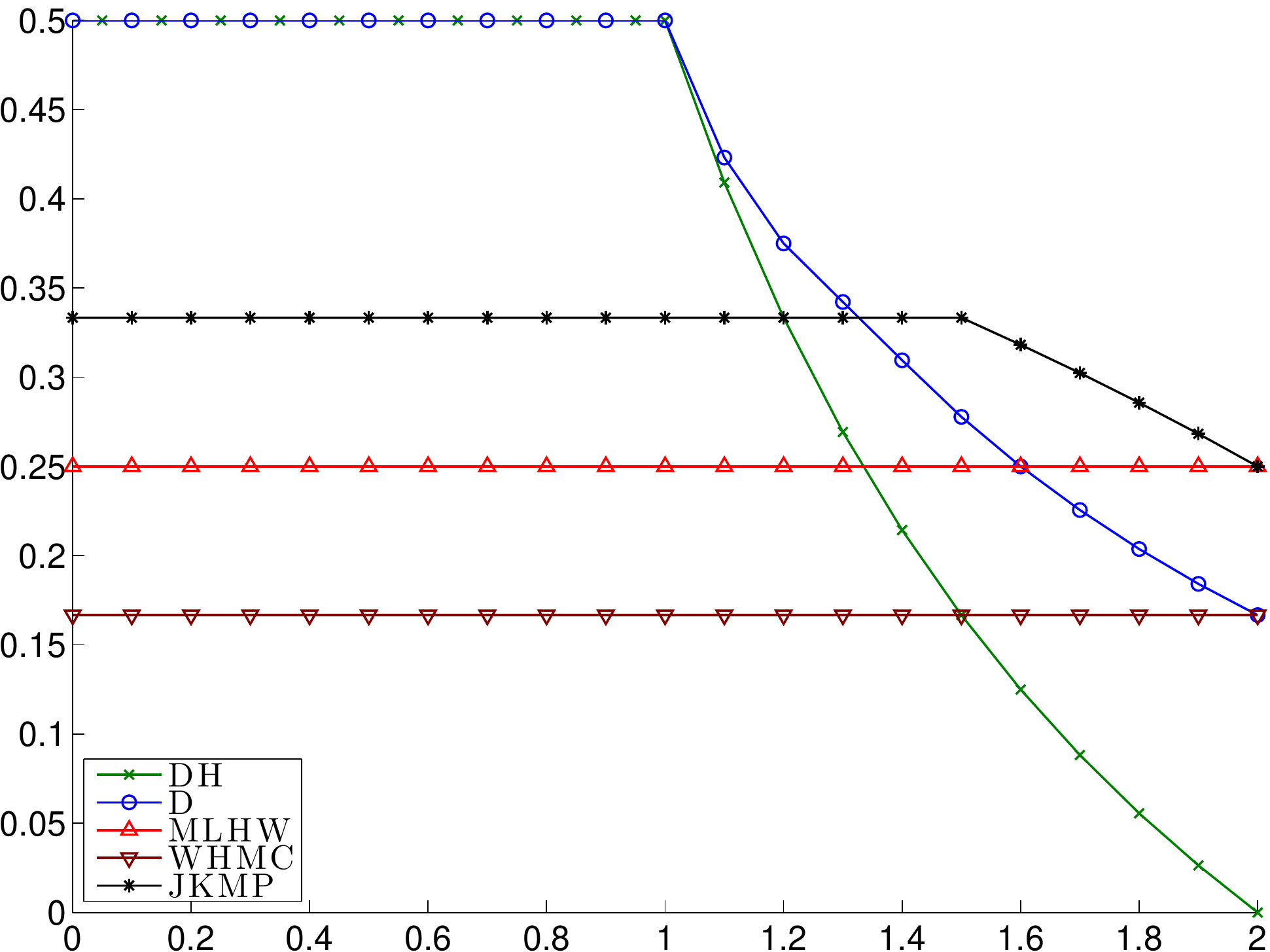}
\caption{\label{rates} Order of convergence with respect to the Blumenthal--Getoor index. }
\end{figure}

The approach in Dereich and  Heidenreich \cite{DH11} is to approximate the underlying L\'evy process via a jump diffusion process, to remove all the small jumps and to {replace them by an additional linear component}. While this method works optimally for processes with jump component of bounded variation, i.e. Blumenthal-Getoor index $\rho \le 1$, the convergence rate degenerates rapidly when $\rho$ gets larger and it becomes arbitrarily bad as $ \rho \to 2$. See Figure~\ref{rates} for a plot of the rate of convergence of this method (denoted \textbf{DH} in the figure). As shown in Dereich \cite{D11}, a better performance is possible when a Brownian correction is added in addition to the linear trend to better take account of the discarded small jumps. This is the curve denoted by \textbf{D} in Figure~\ref{rates}. It converges even for $\rho=2$, but the convergence rate in this case decreases to $\nu^{-\frac{1}{6}}$.

Our Wiener-Hopf based methods, \textbf{WHMC} (single--level) and \textbf{MLWH} (multilevel), are indifferent to the type of jump structure and thus entirely independent of the Blumenthal-Getoor index $\rho$. The price to pay  for this robustness is that, in the bounded variation regime (\ie for L\'evy processes without Gaussian component and with $\rho < 1$), the multilevel version is not quite optimal; recall that for a very general class of L\'evy processes $X$ and payoff functions $F$ the optimal rate of convergence is $\nu^{-\frac{1}{2}}$ (cf. Creutzig et al. \cite{CDMR09} and Dereich \cite{D11} for comments to that effect). Indeed, for $\rho<1$ and $\sigma=0$, the multilevel algorithm performance does not match the rate achieved by the methods of Dereich, but it coincides with the rate reported in Jacod et al. \cite{JKMP05} (depicted as \textbf{JKMP} in Figure \ref{rates}). The performance of the single level algorithm in the bounded variation regime improves to the rate of \textbf{MLWH} for the general case. We do not depict these lines in Figure \ref{rates}. Let us instead focus on the regime where the Blumenthal-Getoor index is $\rho>1$, \ie the L\'evy process has unbounded variation paths. In this case our single--level algorithm \textbf{WHMC} performs better than the method in Dereich and Heidenreich \cite{DH11} for $\rho>4/3$, but it is not able to beat the improved method in Dereich \cite{D11}. %replacement of the small jumps by a Gaussian component adopted in Dereich \cite{D11}.
The multilevel version \textbf{MLWH}, on the other hand, does outperform Dereich \cite{D11} for $\rho>8/5$ and seems to be at least as good as the best currently available method for $\rho =2$. Let us reiterate though that the class of functions $F$ we have so far analysed is particular in that they depend only on the terminal value $X_t$ and the maximum $\overline{X}_t$ of $X$, whereas in Dereich \cite{D11} and Dereich and Heidenreich \cite{DH11} more general functionals of the entire path of $X$ are allowed.
Nevertheless, in the region where the Blumenthal-Getoor index is close to $2$ our multilevel version \textbf{MLWH} is better or equal to any other existing algorithm. This becomes important when we do not have any \emph{a priori} information of the L\'evy process to be used in a numerical application, and hence a robust method with respect to the jump structure is crucial. For example, one often faces the problem of fitting a parametric family of L\'evy processes to given data with a range of parameters allowing unbounded variation jump component. In such cases the method proposed by Dereich and  Heidenreich \cite{DH11} may not convergence in any reasonable time since the numerical performance of the method is restricted by the worst convergence rate in Figure \ref{rates}.

As mentioned above, we believe that our approach can also be adapted to the more general scenario considered in those publications, involving SDEs and path dependent functionals $F$, and still retain the same rates of convergence. We leave this issue to a further piece of work. In the context of more special functionals, it is also worth recalling the results in Jacod et al.~\cite{JKMP05}. Our multilevel algorithm only matches the performance of their method for $\rho=2$. % and Dereich \cite{D11} for $\rho > \frac{4}{3}$. 
However, the assumptions in \cite{JKMP05} are more restrictive than our assumptions and those in Dereich et al. \cite{DH11}, since only functionals of the terminal value $X_t$ are allowed. Moreover, their functions are required to be four times continuously differentiable, highlighting some of the drawbacks and advantages when considering weak instead of strong estimates (cf.~Remark \ref{weak_vs_strong_error}).  

The next section is devoted to a description of another Wiener-Hopf-based algorithm, which achieves optimal convergence rates up to a logarithmic term for all values of the Blumenthal--Getoor index. Despite the fact that we are able to derive  theoretical rates of convergence, the algorithm requires the ability to sample from a triplet of dependent distributions coming out of the Wiener-Hopf factorisation, for which currently no known significant examples are known.

\section{The random time $\mathbf{T}(n,t)$}\label{randomisation_time}

A careful look at Section \ref{sect.WHMC.method} reveals that there is nothing special about the random time $\mathbf{g}(n,n/t)$ in the description of the WHMC method of Theorem \ref{thm_main} other than it being a sum of independent and exponentially distributed random variables converging almost surely  to $t$. Therefore one may investigate more efficient ways to approximate $t$ by a sum of exponentials. We now present an alternative random time that %will turn out to 
exhibits faster convergence rates.

Let $N:=\{N_{s}:s\geq0\}$ be a Poisson process independent of $X$ with intensity $t/n$, where $t>0$ and $n\in\mathbb{N}$, associated with the Wiener-Hopf random walk (\ref{1}). Let $T^N_{j}$ be the time of the $j$-th arrival in $N$, that is to say $T^N_j = \inf\{s>0 : N_s =j\}$ and define  $\mathbf{T}(n,t)=T^N_{N_{t}+1}$. Then, it is well known and intuitively clear from the lack of memory property of the exponentially distributed inter-arrival times of $N$, that $\mathbf{T}(n,t)-t$ follows an exponential distribution with mean $t/n$. Moreover, 
$\lim_{n\uparrow\infty}\mathbf{T}(n,t) = t$ almost surely, and hence  $(X_{\mathbf{T}(n,t)} , \overline{X}_{\mathbf{T}(n,t)} ) $ also converges to $(X_{t}, \overline{X}_{t})$ almost surely. Therefore the same arguments as those found in Theorem~\ref{thm_main} suggest that one may try to approximate $(X_{\mathbf{T}(n,t)}, \overline{X}_{\mathbf{T}(n,t)})$ by $(V(N_{t}+1,n/t), J(N_{t}+1,n/t))$. Unfortunately, current knowledge of the Wiener-Hopf factorisation only gives us  information about the pair
$(\overline{X}_{\mathbf{e}(q)},X_{\mathbf{e}(q)}-\overline{X}_{\mathbf{e}(q)})$, rather than the triplet $(\overline{X}_{\mathbf{e}(q)},X_{\mathbf{e}(q)}-\overline{X}_{\mathbf{e}(q)},\mathbf{e}(q))$. Therefore we are not able to perform a simultaneous exact sample from $(V(N_{t}+1,n/t), J(N_{t}+1,n/t),N_{t}+1)$, hence the random time $\mathbf{T}(n,t)$ cannot be used in a practical implementation with the theory developed in Kuznetsov et al. \cite{KKPvS}. 

The time horizon $\mathbf{T}(n,t)$ leads to a slightly different construction or modification of the WHMC estimator  described in Theorem \ref{thm_main} and eventually a different construction of the MLWH estimator and its associated algorithm, but for the sake of simplicity 
%we will refer to each construction or algorithm by the time horizon which describes it. The 
the Monte-Carlo estimate of $\mathbb{E}(F(\overline{X}_t, X_t))$ is still denoted by  $\widehat{F}_{\mathrm{MC}}^{n,M}$ and written as in (\ref{WH-MC}). However, for any $i\in\mathbb{N}$,  we now understand $F^{n,(i)}$ as the $i$-th draw of the random variable
\begin{equation}\label{general_estimator}
F^{n} =F(V(N_{t}+1, n/t), J(N_{t}+1, n/t))\ .
\end{equation}
\begin{remark}\em
A third time horizon which constitutes only a minor variation of $\mathbf{T}(n,t)$, is to consider the random time $\widetilde{\mathbf{T}}(n,t):=T^N_{N_{t}}$, in which case $t-\widetilde{\mathbf{T}}(n,t)$ has the same distribution as $\mathbf{e}(n/t)\wedge t$. It is an easy exercise to see that all results applying to $\mathbf{T}(n,t)$ will hold also for $\widetilde{\mathbf{T}}(n,t)$ and therefore we avoid any further comments on $\widetilde{\mathbf{T}}(n,t)$.
\end{remark}

\subsection{Thinning and convergence theorems}

For the time horizon $\mathbf{T}(n,t)$ the thinning procedure follows analogously to the description in Section \ref{refinement} and equations (\ref{refinement_1}) and (\ref{refinement_2}) remain true. The issue of  consistency is not entirely straight forward however. Following the notation in Section \ref{refinement}, let $\ell$ be one of the levels in the multilevel algorithm and denote by $\check{\mathbf{T}}(n_{\ell-1},t)$ the time horizon produced by thinning the Poisson process $N^{\ell}$. We now show that $\check{\mathbf{T}}(n_{\ell-1},t)$ has the same distribution as ${\mathbf{T}}(n_{\ell-1},t)$. It is clear from the algorithm that
\begin{equation}\label{thinning_T}
{\mathbf{T}}(n_{\ell},t)\leq\check{\mathbf{T}}(n_{\ell-1},t)\ .
\end{equation}
Indeed, due to the coin tossing either the two time horizons in (\ref{thinning_T}) are equal or differ in a geometric sum of exponentials. Since $\mathbf{T}(n_{\ell},t)-t\stackrel{d}{=}\mathbf{e}(n_{\ell}/t)$ we note that
\begin{equation}\label{consistency_T}
\check{\mathbf{T}}(n_{\ell-1},t)-t=
\mathbf{T}(n_{\ell},t)-t+\sum_{i=1}^{\kappa}\mathbf{e}_{i}(n_{\ell}/t)
\stackrel{d}{=}\sum_{i=0}^{\kappa}\mathbf{e}_{i}(n_{\ell}/t)
\stackrel{d}{=}\mathbf{e}(n_{\ell-1}/t)\stackrel{d}{=}{\mathbf{T}}(n_{\ell-1},t)-t\ ,
\end{equation}
where $\kappa$ is an independent geometric distribution on $\{0,1,2,\cdots\}$ with parameter $1/2$.
We then conclude that, as in Section \ref{refinement}, $\check{F}^{n_{\ell-1}}$ and ${F}^{n_{\ell-1}}$ have the same law.

Equations (\ref{thinning_T}) and (\ref{consistency_T}) suggest an interesting simplification. Given the pair $(X_{\mathbf{T}(n_\ell,t)} , \overline{X}_{\mathbf{T}(n_\ell,t)})$, rather than extending the random grid on level $\ell$ and using the thinning algorithm together with (\ref{refinement_1}) and (\ref{refinement_2}) to produce a coarser sample, we only have to flip a coin once. If it is a head, then $(X_{\mathbf{T}(n_{\ell-1},t)}^{(i)} , \overline{X}^{(i)}_{\mathbf{T}(n_{\ell-1},t)} )=(X^{(i)}_{\mathbf{T}(n_\ell,t)} , \overline{X}^{(i)}_{\mathbf{T}(n_\ell,t)} )$ and thus the entire sample $F^{n_\ell,(i)}-F^{n_{\ell-1},(i)} = 0$. On the other hand, if it is a tail, we just need to simulate one independent draw $(S^{(i)}_{n_{\ell-1}/t},I^{(i)}_{n_{\ell-1}/t})$ of $(\overline{X}_{\mathbf{e}(n_{\ell-1}/t)},\underline{X}_{\mathbf{e}(n_{\ell-1}/t)})$ and set 
\begin{eqnarray}
\overline{X}^{(i)}_{\mathbf{T}(n_{\ell-1},t)}&=&\overline{X}^{(i)}_{\mathbf{T}(n_{\ell},t)} \, \vee \, (X^{(i)}_{\mathbf{T}(n_{\ell},t)}+S^{(i)}_{n_{\ell-1}/t}) \;,\label{sup_ext_thin}\\
X^{(i)}_{\mathbf{T}(n_{\ell-1},t)}&=&X^{(i)}_{\mathbf{T}(n_{\ell},t)}+S^{(i)}_{n_{\ell-1}/t}+I^{(i)}_{n_{\ell-1}/t}\label{term_ext_thin}\ .
\end{eqnarray}
Observe that equations (\ref{sup_ext_thin}) and (\ref{term_ext_thin}) significantly simplify the multilevel algorithm for the time horizon $\mathbf{T}(n,t)$.  %much faster than using equations (\ref{refinement_1}) and (\ref{refinement_2}).

%We are now ready to derive the convergence results for the WHMC and MLMC algorithms adapted to the time horizon $\mathbf{T}(n,t)$.

\begin{prop}\label{mean_square_error_Poisson}
Let $t>0$ and suppose $X$ satisfies Assumption (A1). Then, for any $n \in \mathbb{N}$, %we have
\begin{enumerate}
\item[(i)] $\E[(X_{\mathbf{T}(n,t)}-X_{t})^{2}]\eqsim n^{-1}$ \ and \ $\E[|X_{\mathbf{T}(n,t)}-X_{t}|]\lesssim n^{-1/2}$,
\item[(ii)] $\E[(\overline{X}_{\mathbf{T}(n,t)}-\overline{X}_{t})^{2}]\lesssim n^{-1}$ \ and \ $\E[|\overline{X}_{\mathbf{T}(n,t)}-\overline{X}_{t}|]\lesssim n^{-1/2}$.
\end{enumerate}
\end{prop}
\proof As in the proof of Proposition \ref{mean_square_error_gamma}, it is enough to prove the second moment results in (i) and (ii), since the results for the first moments then follow from Jensen's inequality.  Thanks to  the lack of memory property, it is well known that $\mathbf{T}(n,t)-t$ follows an exponential distribution with mean $t/n$.
It follows that, for all $n\in\mathbb{N}$,
\begin{equation}
\E[(\mathbf{T}(n,t)-t)^{2}]\;=\;\frac{2t^{2}}{n^{2}} \quad \text{and} \quad %\label{poisson_1}\\
\E[|\mathbf{T}(n,t)-t|]\;=\;\E[\mathbf{T}(n,t)-t]=\frac{t}{n}\label{poisson_1_T}\ .
\end{equation}
 The desired result now follows by appealing to Lemma \ref{mean_square_error_general}. 
\qed 

As for Theorem \ref{g-error}, the next theorem follows immediately from Theorems \ref{thm_WHMC} and \ref{Rob_thm}, Proposition \ref{mean_square_error_Poisson}, and the fact that $\gamma = 1$ (cf.~Assumption (A3)). Note that in this case the cost of producing a single path $F^{n_\ell}$ is random, as the Wiener-Hopf walk will have a random number of iterations. It is nevertheless clear that $\gamma=1$ in Theorem \ref{thm_WHMC}, but also in Theorem \ref{Rob_thm}, since we will have to produce at most one additional exponential time step with rate $n_{\ell-1}$ in order to construct the coarser sample (as shown in equation (\ref{consistency_T})).

\begin{teo}\label{T-error}
Suppose that assumptions (A1)--(A3) are satisfied and that the random time horizon is $\mathbf{T}(n,t)$. Then, for any $\nu \in \mathbb{N}$ and under the constraint that the expected value of the total computational cost is $\mathcal{O}(\nu)$ operations, the root mean square error for the WHMC and the MLWH method can be bounded respectively by 
\begin{equation*}%\label{rate_T}
e(\widehat{F}_{\mathrm{MC}}^{n,M})\lesssim \nu^{-\frac{1}{4}} \qquad \text{and} \qquad e(\widehat{F}_{\mathrm{ML}}^{\mathcal{M}(n_0,L)})\lesssim \nu^{-\frac{1}{2}} \log^2\nu.
\end{equation*}
%(with suitable choices for $n$, $M$, $L$ and $\mathcal{M}(n_0,L)$, as stated in Theorems \ref{thm_WHMC} and \ref{Rob_thm}).
\end{teo}

Note that the multilevel algorithm with time horizon $\mathbf{T}(n,t)$ achieves optimal convergence rate up to a logarithm factor. Recalling the discussion in Section \ref{g_time} and Figure \ref{rates}, the above methodology would offer optimal convergence throughout the entire range of the Blumenthal--Getoor index, thereby beating any of the methods considered there. %Even if the algorithm is not implementable in general there are some trivial examples where the triplet $(\overline{X}_{\mathbf{e}(q)},X_{\mathbf{e}(q)}-\overline{X}_{\mathbf{e}(q)},\mathbf{e}(q))$ is available such as in the Brownian motion case.

\section{Implementation and Numerical Results}\label{Numerics}

The aim of this section is to provide an example of the application of the WHMC and the new MLWH method to demonstrate the feasibility of both algorithms and to confirm the theoretical results of the preceding sections. 

We will study $X_t$ and $\overline{X}_t$, as well as the price of a barrier option assumed to belong to a market driven by a member of the $\beta$-class. We use this class because there is an explicit form of the Wiener-Hopf factorisation and because, through appropriate parameter choices, it contains examples of processes with paths of bounded and unbounded variation. Moreover % and because it can exhibit the entire variety of behaviours of L\'evy processes through a change of parameters. In fact,  
Kuznetsov et al. \cite[Section 4]{KKPvS} argue that a large class of L\'evy processes can be approximated by a member of the $\beta$-class plus a compound Poisson process. It is beyond  the scope of this paper to discuss the $\beta$-class in a financial framework, but we refer the reader to Ferreiro-Castilla and Schoutens \cite{FS11} and Schoutens and van Damme \cite{SD10} where the $\beta$-class is studied in relation with other popular models in finance and where a comprehensive numerical analysis of the computation of the extrema for this family is done. We will simply introduce the $\beta$-class as a ten parameter L\'evy process with triplet $(a,\sigma, \Pi)$, such that the L\'evy measure is absolutely continuous with density
\begin{equation}\label{beta_density}
\Pi({\rm d}x)= c_1\frac{{\mathrm e}^{-a_1 b_1 x}}{(1- {\mathrm e}^{-b_1 x})^{\lambda_1}}\mathbf{1}_{\{x>0\}}{\rm d}x
+ c_2 \frac{{\mathrm e}^{a_2 b_2 x}}{(1- {\mathrm e}^{b_2 x})^{\lambda_2}}\mathbf{1}_{\{x<0\}}{\rm d}x\ ,
\end{equation}
and $a_{i}>0$, $b_{i}>0$, $c_{i}\geq0$ and $\lambda_{i}\in(0,3)$, for $i=1,2$. The expressions for the distribution of the extrema can be found in Kuznetsov \cite[Section 4]{Kuz} and are given by the following infinite product representation
\begin{equation*}
\mathbb{E}[e^{iz\overline{X}_{\mathbf{e}(q)}}]=\frac{1}{1+\frac{iz}{\zeta_{0}^{-}(q)}}\prod_{n\leq-1}\frac{1+\frac{iz}{b_{1}(n+1-a_{1})}}{1+\frac{iz}{\zeta_{n}(q)}}
\qquad\text{and}\qquad
\mathbb{E}[e^{iz\underline{X}_{\mathbf{e}(q)}}]=\frac{1}{1+\frac{iz}{\zeta_{0}^{+}(q)}}\prod_{n\geq1}\frac{1+\frac{iz}{b_{2}(n-1+a_{2})}}{1+\frac{iz}{\zeta_{n}(q)}}\ ,
\end{equation*}
where $\zeta_{0}^{-}(q)$ and $\{\zeta_{n}(q)\}_{n\leq-1}$ correspond to the negative zeros and $\zeta_{0}^{+}(q)$ and $\{\zeta_{n}(q)\}_{n\geq1}$ correspond to the positive zeros of $\Psi(-iz)+q=0$.
The numerical simulation of $\overline{X}_{\mathbf{e}_q}$ and $\underline{X}_{\mathbf{e}_q}$ requires truncating the infinite products above, but it can be shown using \cite[Theorem 10]{Kuz} and \cite[Section 3.1]{FS11} that the cost to obtain a sample via this process (to an arbitrary degree of accuracy) is independent of the value of $q$, thus verifying Assumption (A3). In fact this claim holds even for the more general family of processes introduced in Kuznetsov et al. \cite{KKPvS} known as meromorphic L\'evy processes. From (\ref{beta_density}), one sees that processes belonging to the $\beta$-class have L\'evy measures which decay exponentially and therefore have finite moments of all orders, in particular they verify Assumption (A1) for the entire range of parameters. 

For the numerical results below, we choose $a=c_{i}=a_{i}=b_{i}=1$ and set $\lambda_{2}=\lambda_{1}=\lambda$. We then study the behaviour and the robustness of our methods for various choices of $\sigma$ and $\lambda$. For instance, if we set $\lambda=0.5$ then we have a process of finite activity and thus with Blumenthal--Getoor index $\rho=0$. By letting $\lambda=1.5$, we have an infinite activity process with jump component of bounded variation (i.e. $\rho \in (0,1)$), and finally by setting $\lambda=2.5$, we have an infinite activity process with index $\rho>1$. We are free to chose any value of the Gaussian coefficient $\sigma$.  When $\sigma=0$ only the case  with $\lambda=2.5$ will produce a process with paths of unbounded variation. When $\sigma^{2}>0$ all processes have paths of unbounded variation.

The pay-off function in our case study relates to the pricing of barrier options, \ie we are interested in computing 
\begin{equation}\label{output_functional}
\E[(K-\exp(x_{0}+X_{t}))^{+}\mathbf{1}_{\{x_{0}+\overline{X}_{t}>b\}}]
\end{equation}
for some $K,b>0$ and $x_0 \in \mathbb{R}$. Barrier options play an important role in the computation of exotic options and are the building blocks for more complex derivatives such as Credit Default Swaps, see for instance Schoutens and Cariboni \cite{SC}. Even though the pay-off function in \eqref{output_functional} is not Lipschitz continuous, and thus does not satisfy Assumption (A2), we will see below that the theoretically predicted performance is nevertheless achieved.

The programs were implemented in C++, in long double precision and executed on a 64-bit Intel(R) Xeon(R) CPU at 2.66GHz running 32-bit Linux 2.6.32-41-generic-pae using the GNU gcc compiler (version 4.8). Although currently running in a single threaded environment, (multilevel) Monte Carlo algorithms are of course ideally posed to be parallelised, as individual samples and simulations at different levels are independent of each other.

\subsection{Numerical results}\label{experiments}

\begin{figure}[t!]
\center
\begin{minipage}[t]{0.47\textwidth}
\begin{flushright}
\includegraphics[height=6.8cm,width=8cm]{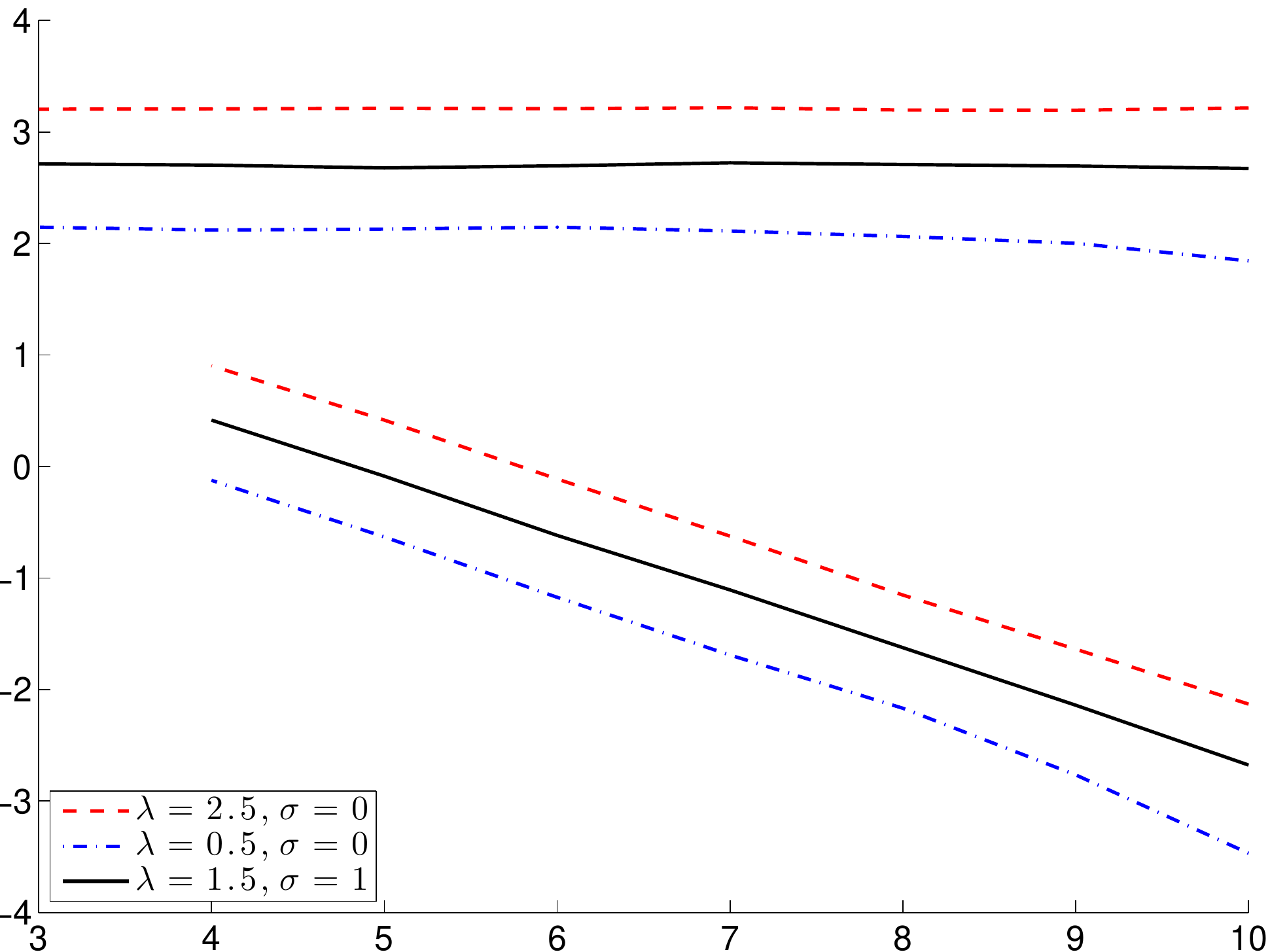}
\end{flushright}
\end{minipage}
\hskip0.03\textwidth
\begin{minipage}[t]{0.47\textwidth}
\begin{center}
\includegraphics[height=6.8cm,width=8cm]{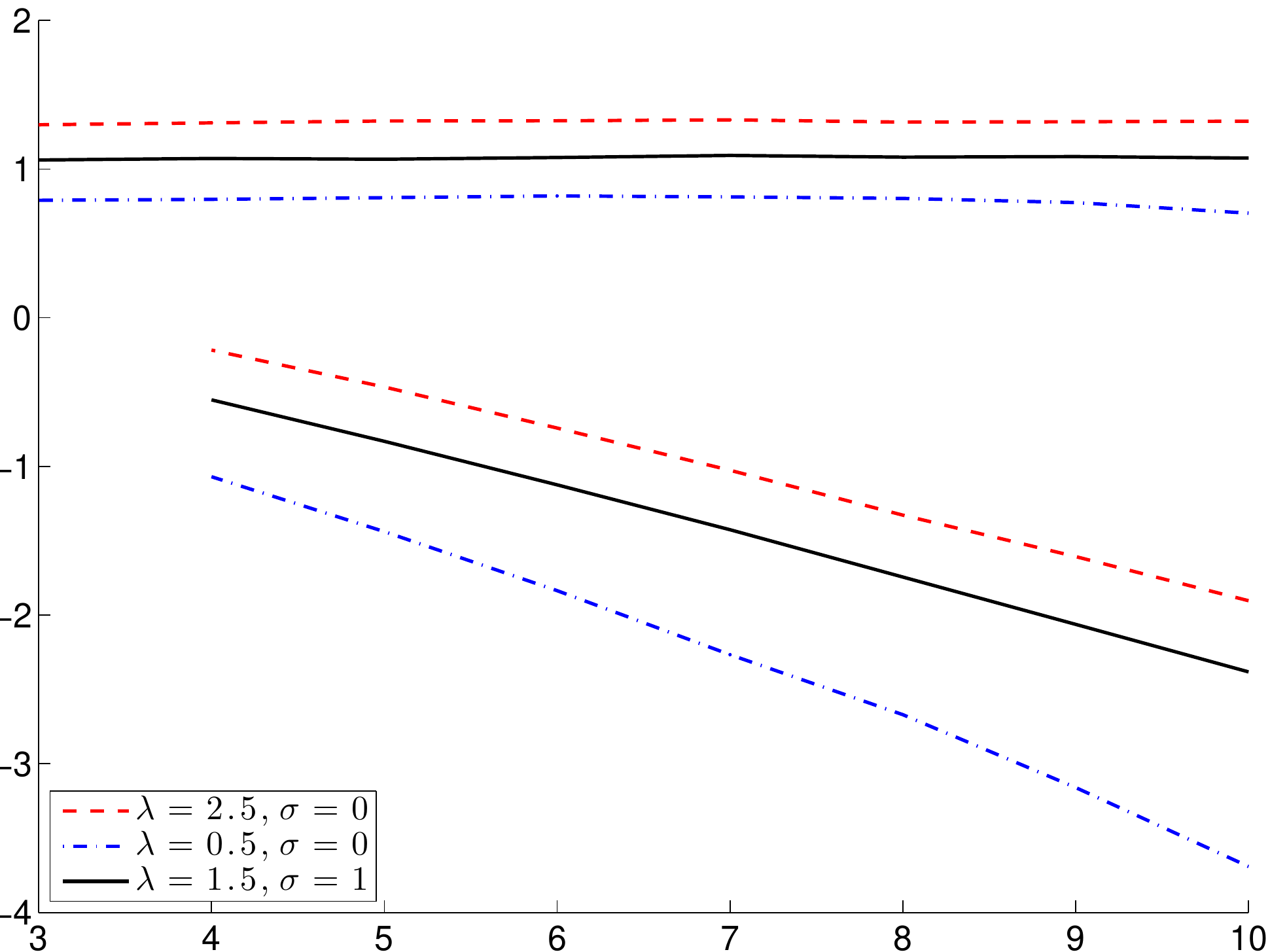}
\end{center}
\end{minipage}
\caption{\label{plot1} $\V(X_{\mathbf{g}(n_\ell,n_\ell/t)}-X_{\mathbf{g}(n_{\ell-1},n_{\ell-1}/t)})$ and $\V(X_{\mathbf{g}(n_\ell,n_\ell/t)})$ (left figure), as well as $\E[|X_{\mathbf{g}(n_\ell,n_\ell/t)}-X_{\mathbf{g}(n_{\ell-1},n_{\ell-1}/t)}|]$ and $\E[|X_{\mathbf{g}(n_\ell,n_\ell/t)}|]$ (right figure) plotted against $n_\ell$ (using $\log_2$--scales).}
\end{figure}

\begin{figure}[t!]
\center
\begin{minipage}[t]{0.47\textwidth}
\begin{flushright}
\includegraphics[height=6.8cm,width=8cm]{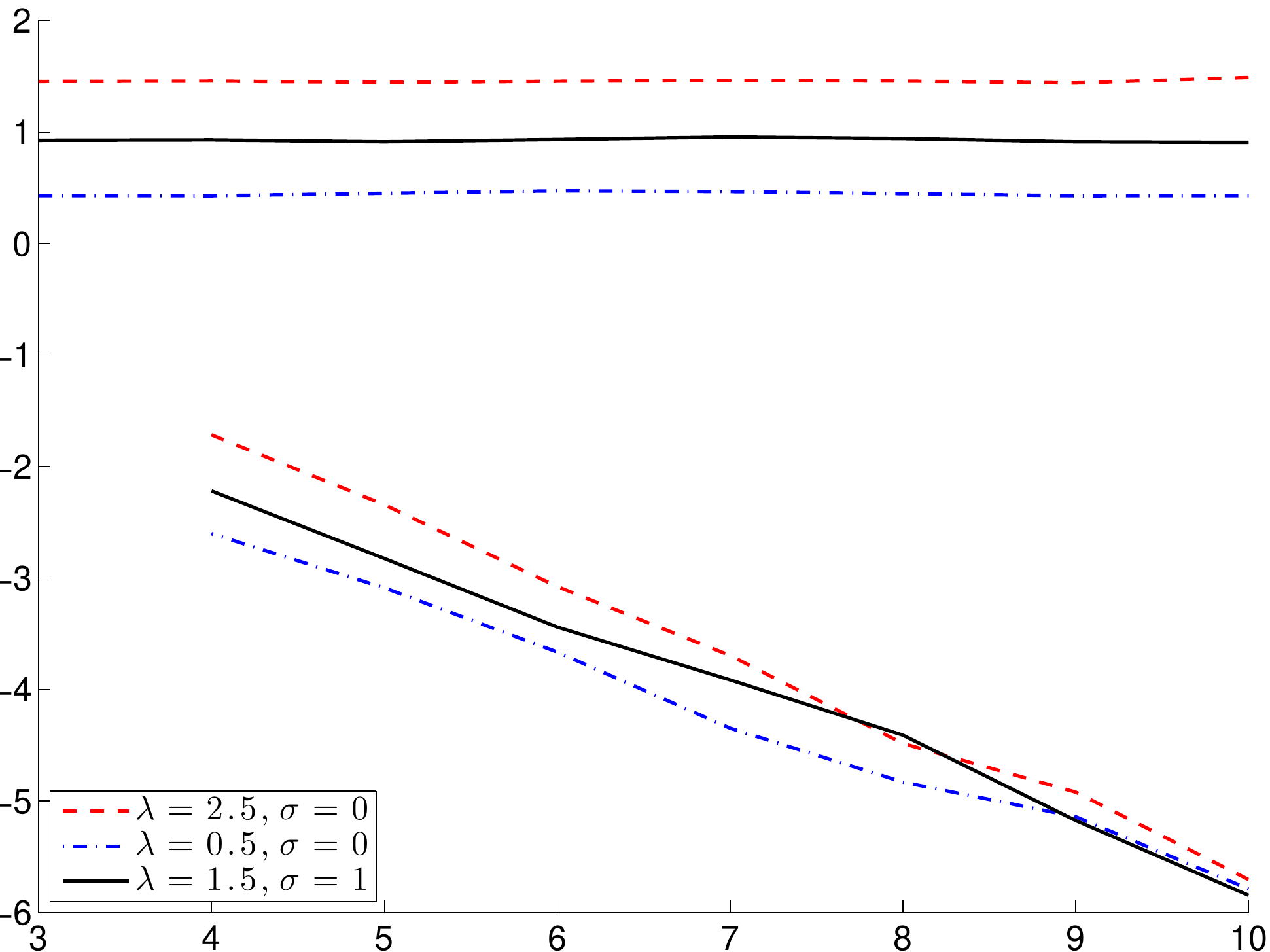}
\end{flushright}
\end{minipage}
\hskip0.03\textwidth
\begin{minipage}[t]{0.47\textwidth}
\begin{flushright}
\includegraphics[height=6.8cm,width=8cm]{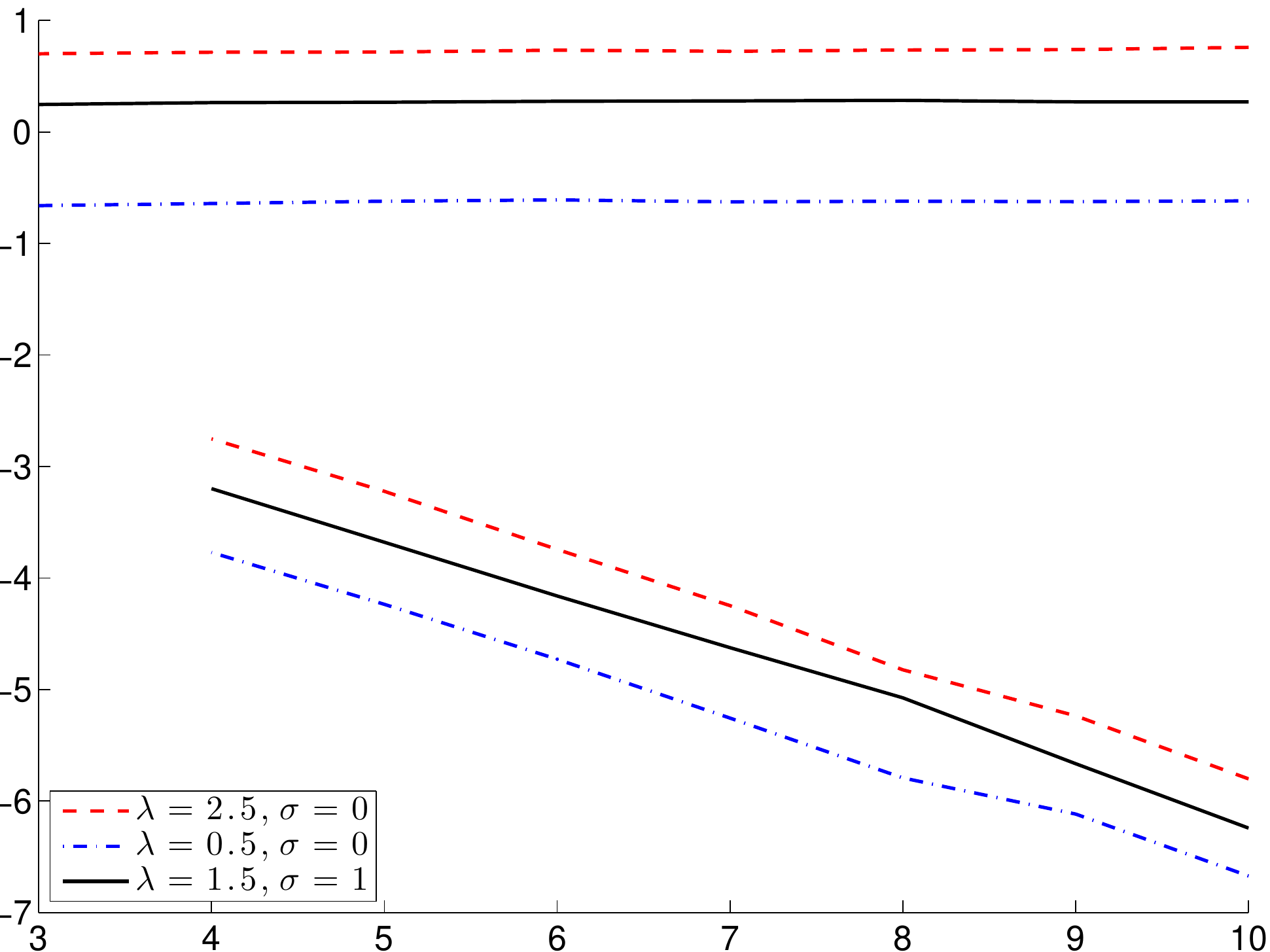}
\end{flushright}
\end{minipage}
\caption{\label{plot2} $\V(\overline{X}_{\mathbf{g}(n_\ell,n_\ell/t)}-\overline{X}_{\mathbf{g}(n_{\ell-1},n_{\ell-1}/t)})$ and $\V(\overline{X}_{\mathbf{g}(n_\ell,n_\ell/t)})$ (left figure), as well as $\E[|\overline{X}_{\mathbf{g}(n_\ell,n_\ell/t)}-\overline{X}_{\mathbf{g}(n_{\ell-1},n_{\ell-1}/t)}|]$ and $\E[|\overline{X}_{\mathbf{g}(n_\ell,n_\ell/t)}|]$ (right figure) plotted against $n_\ell$ (using $\log_2$--scales).}
\end{figure}

 We start in Figures~\ref{plot1} and \ref{plot2} by plotting, for various pairs of $\lambda$ and $\sigma$ and for $n_\ell = 2 n_{\ell-1}$,  the mean and the variance of 
$(X_{\mathbf{g}(n_\ell,n_\ell/t)}-X_{\mathbf{g}(n_{\ell-1},n_{\ell-1}/t)})$ and $(\overline{X}_{\mathbf{g}(n_\ell,n_\ell/t)}-\overline{X}_{\mathbf{g}(n_{\ell-1},n_{\ell-1}/t)})$. For reference we also include mean and variance of $X_{\mathbf{g}(n_\ell,n_\ell/t)}$ and $\overline{X}_{\mathbf{g}(n_\ell,n_\ell/t)}$. We see from the left plots in those figures that the convergence of $\V(X_{\mathbf{g}(n_\ell,n_\ell/t)}-X_{\mathbf{g}(n_{\ell-1},n_{\ell-1}/t)})$ 
%and $\V(\overline{X}_{\mathbf{g}(n_\ell,n_\ell/t)}-\overline{X}_{\mathbf{g}(n_{\ell-1},n_{\ell-1}/t)})$ 
with respect to $n_\ell^{-1/2}$ is linear as predicted in Proposition \ref{mean_square_error_gamma} and slightly better than predicted for $\V(\overline{X}_{\mathbf{g}(n_\ell,n_\ell/t)}-\overline{X}_{\mathbf{g}(n_{\ell-1},n_{\ell-1}/t)})$. % in Figure \ref{plot2}.

Since the rate of the Cauchy sequence $\E[|X_{\mathbf{g}(n_\ell,n_\ell/t)}-X_{\mathbf{g}(n_{\ell-1},n_{\ell-1}/t)}|]$ dictates the rate for $\E[|X_{\mathbf{g}(n_\ell,n_\ell/t)}-X_{t}|]$, we can also verify that the convergence of $\E[|X_{\mathbf{g}(n_\ell,n_\ell/t)}-X_t|]$ and $\E[|\overline{X}_{\mathbf{g}(n_\ell,n_\ell/t)}-\overline{X}_{t}|]$ is no worse than $\mathcal{O}(n_\ell^{-1/4})$. In fact, in the absence of a Gaussian component and in the case bounded variation jump activity, the convergence of $\E[|X_{\mathbf{g}(n_\ell,n_\ell/t)}-X_{\mathbf{g}(n_{\ell-1},n_{\ell-1}/t)}|]$ seems to be even linear in $n_\ell^{-1/2}$, which is not surprising given the piecewise linear nature of this model. More surprisingly though, the bias for the running maximum seems to  converge faster than predicted for any choice of $\lambda$ and $\sigma$.

\begin{figure}[t!]
\center
\begin{minipage}[t]{0.47\textwidth}
\begin{flushright}
\includegraphics[height=6.8cm,width=8cm]{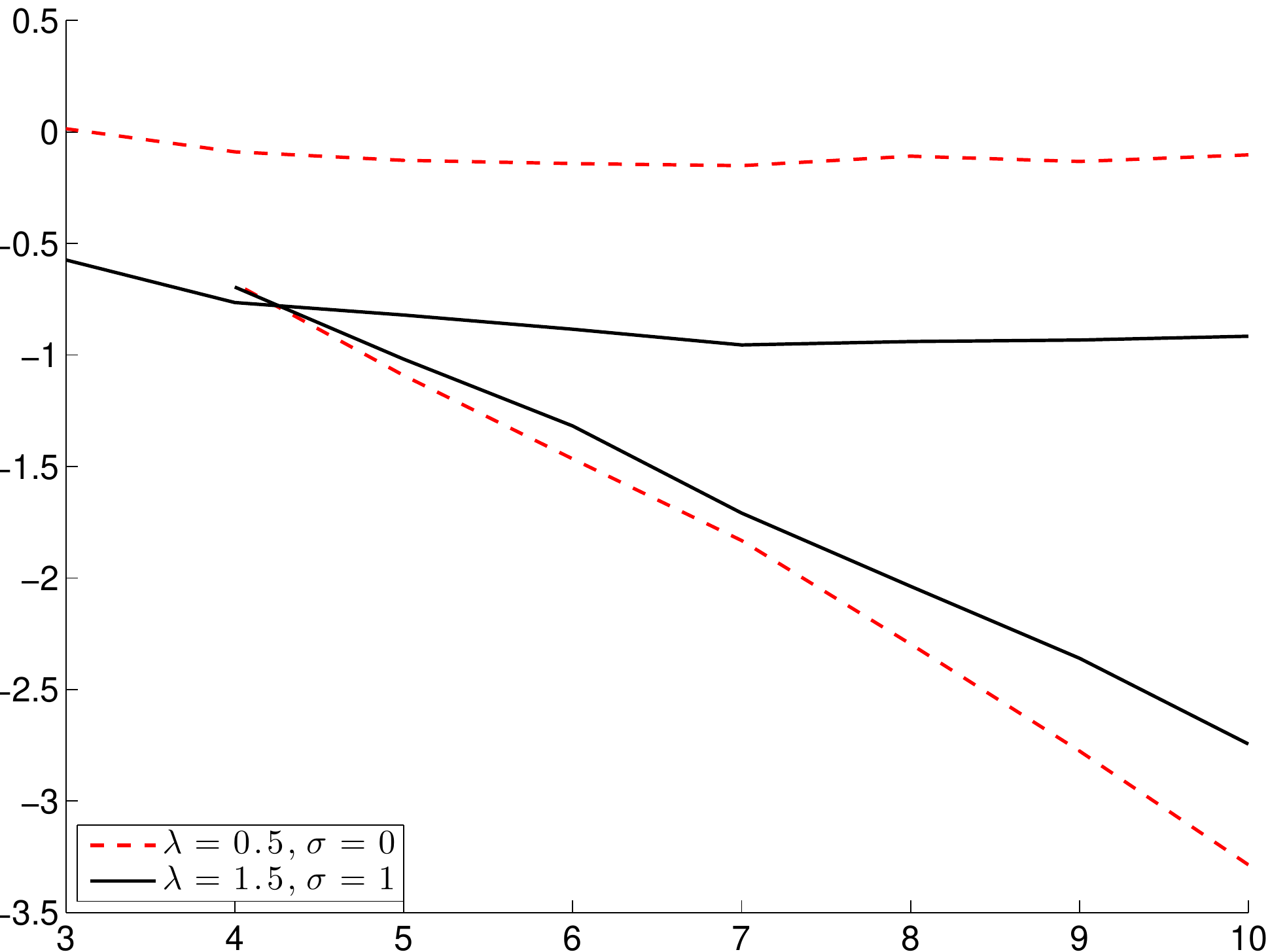}
\end{flushright}
\end{minipage}
\hskip0.03\textwidth
\begin{minipage}[t]{0.47\textwidth}
\begin{center}
\includegraphics[height=6.8cm,width=8cm]{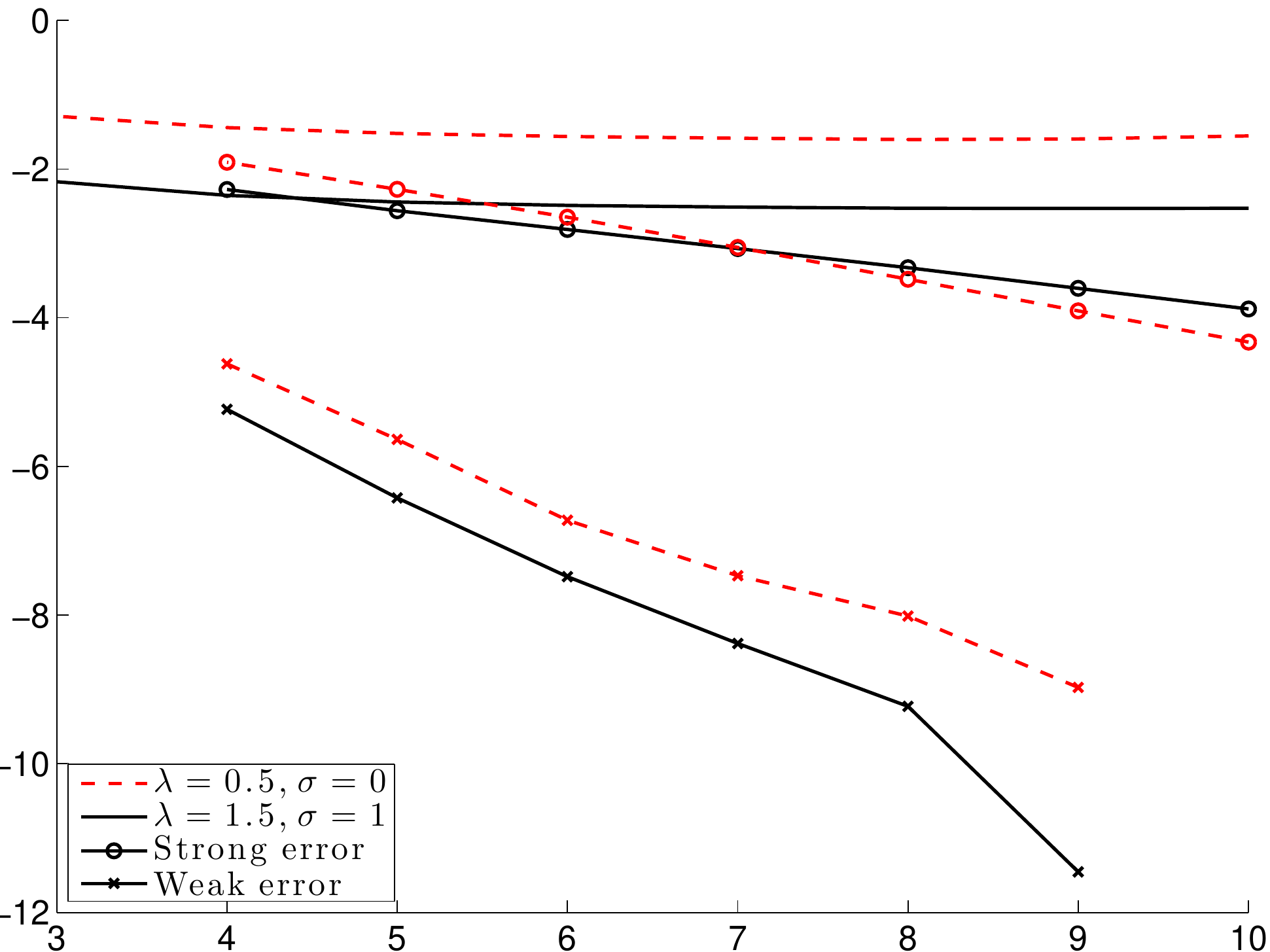}
\end{center}
\end{minipage}
\caption{\label{plot3} $\V(F^{n_\ell}-F^{n_{\ell-1}})$ and $\V(F^{n_\ell})$ (left figure), as well as strong error $\E[|F^{n_\ell}-F^{n_{\ell-1}}|]$,  weak error $|\E[F^{n_\ell}]-\E[F^{n_{\ell-1}}]|$ and $\E(|F^{n_\ell}|)$ (right figure) plotted against $n_\ell$ (using $\log_2$--scales).}
\end{figure}
 
The behaviour for our output functional in (\ref{output_functional}), which violates assumption (A3), is also very close to the conclusions in Theorem \ref{g-error}. As the plots in Figure~\ref{plot3} show, the rate of decay in $\V(F^{n_\ell}-F^{n_{\ell-1}})$ is approximately $0.36$, while Theorem \ref{g-error} claims $\beta=1/2$. Similarly, $\E[|F^{n_\ell}-F^{n_{\ell-1}}|]$ behaves roughly like $\mathcal{O}(n_\ell^{-1/4})$, confirming that $\alpha =1/4$. We also see that $\V(F^{n_\ell})$ remains constant as shown in (\ref{var_cte}). In the right plot of Figure~\ref{plot3} we also include the weak error, \ie $|\E[F^{n_\ell}]-\E[F^{n_{\ell-1}}]|$, which is indeed much smaller than the strong error and decays much faster (cf.~Remark~\ref{weak_vs_strong_error}).

\begin{figure}[t!]
\center
\begin{minipage}[t]{0.47\textwidth}
\begin{flushright}
\includegraphics[height=7cm,width=8cm]{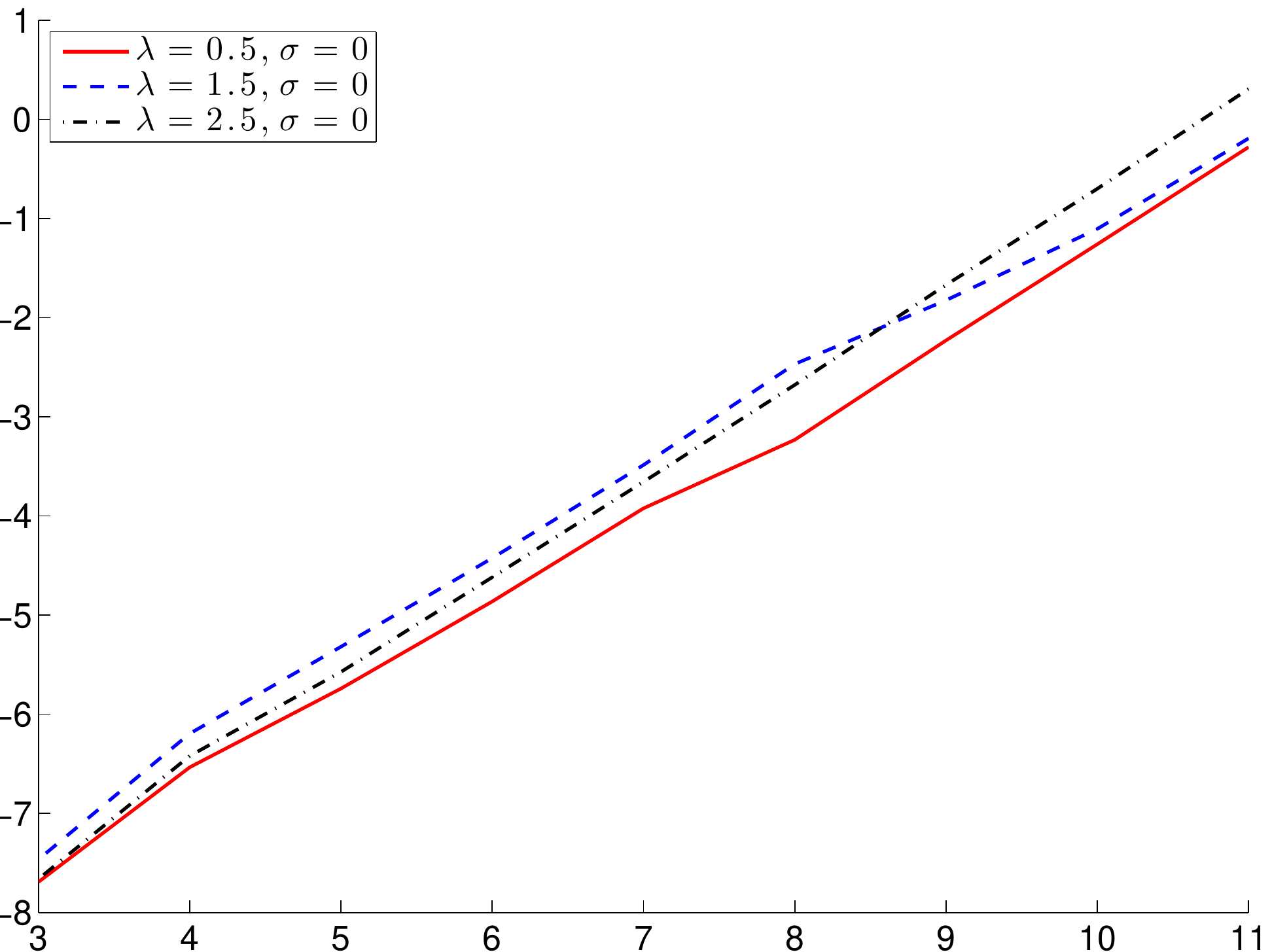}
\end{flushright}
\end{minipage}
\hskip0.03\textwidth
\begin{minipage}[t]{0.47\textwidth}
\begin{center}
\includegraphics[height=7cm,width=8cm]{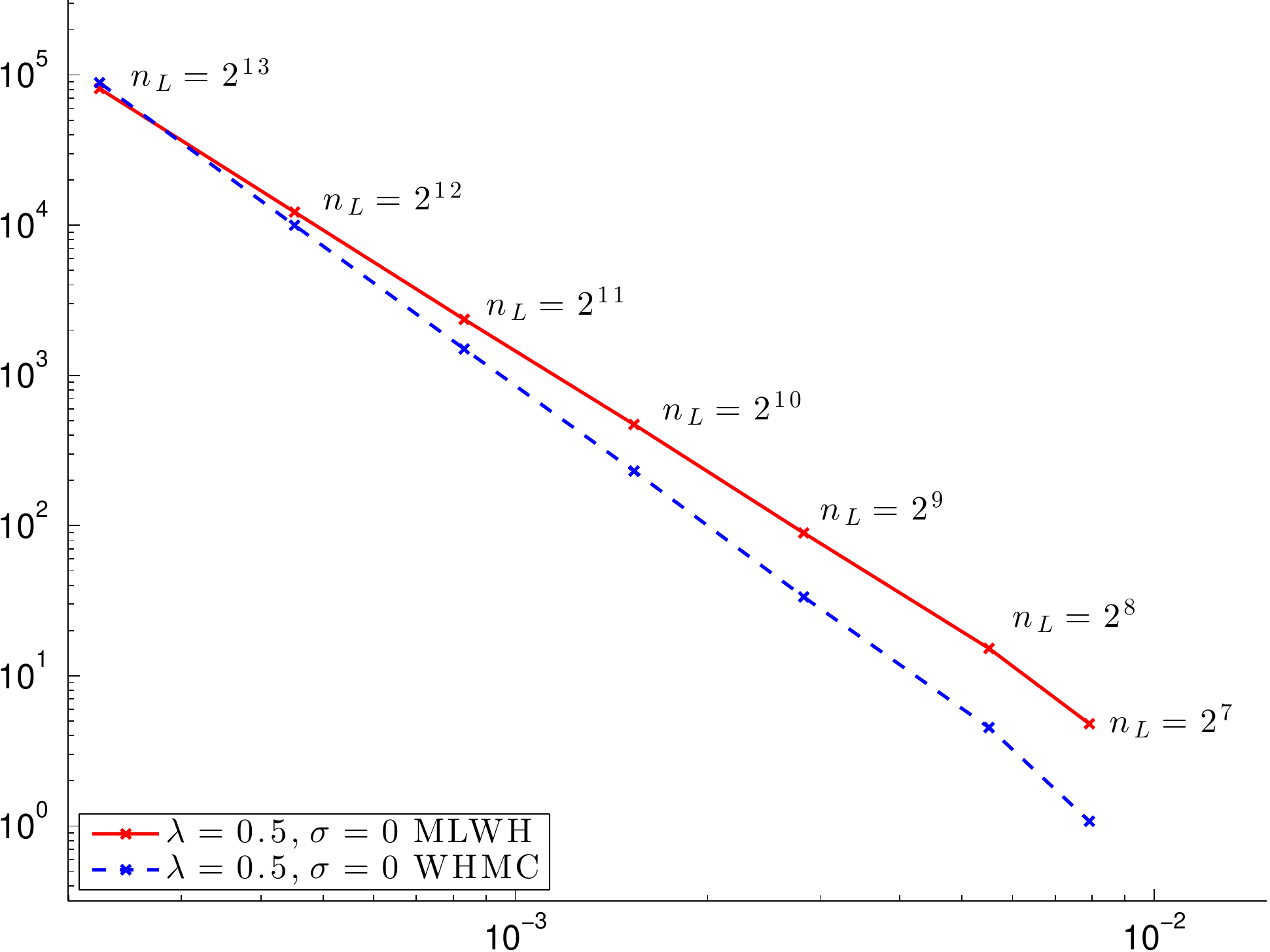}
\end{center}
\end{minipage}
\caption{\label{plot5} Left: Average CPU--time (in milliseconds) to compute one sample $F^{n,(i)}$ plotted against $n = n_L$ ($\log_2$--scales). Right: Total CPU-time (in seconds) for the WHMC and MLWH methods plotted against the estimated root mean square error $e$ ($\log_{10}$--scales).}
\end{figure}

Finally, in Figure \ref{plot5} we look at actual CPU times in running our WHMC and MLWH algorithms. In the left inset of  Figure \ref{plot5} we study the  cost of computing one single sample $F^{n,(i)}$ of $F^{n}$ using the Wiener-Hopf approach, to confirm also numerically Assumption (A3). The cost clearly grows linearly with $n$, independently of any parameters in the L\'evy process. As mentioned in the previous section, this can in fact be proved rigorously.

Our final plot in Figure \ref{plot5} (right) shows actual performance results and compares the single level and the multilevel Wiener-Hopf Monte Carlo algorithms. The $x$-axis describes the estimated root mean square error (RMSE) $e$. The $y$-axis describes the cost $\nu$. In each data point of the plot, the finest level in the MLMC algorithm matches the approximation level in the WHMC. Next to each vertical pair of points we have labeled that value. This ensures that the bias error is the same in both algorithms. The sampling error in the MLMC error is split between levels accordingly to \cite[Equation (12)]{Giles08}. Since the weak error is much smaller than the variance reduction in this particular example, as seen in Figure \ref{plot3}, the multilevel algorithm outperforms the single level version only for relatively small accuracies. However, from Figure~\ref{plot5} one also sees that  the slopes of the two lines clearly play into the hands of the MLMC. We have considered also only one functional so far. The variance curves in Figure \ref{plot2} suggest that for functionals that depend more strongly on the maximum $\overline{X}_t$, the MLWH method would be faster than the single level method already for $n_L = 2^5$ or $n_L = 2^6$. 

Finally, we also have not yet experimented with the factor $s$ at which we coarsen the rates from one level to the next, i.e. $n_\ell = s n_{\ell-1}$. We have only studied the case $s=2$, while in Giles \cite{Giles08} it was shown that the optimal value for standard multilevel Monte Carlo methods for SDEs with Gaussian noise is $s=7$. In any case, the results in Figure \ref{plot5} (right) show that Wiener-Hopf Monte Carlo techniques allow the pricing of options in a market driven by a L\'evy process to an accuracy of $\mathcal{O}(10^{-2})$ in less than a second.\vspace{1ex}

\noindent
{\bf Acknowledgements.} We are grateful to Kees van Schaik and Steffen Dereich for  insightful discussions which lead to significant changes in this paper.

\end{document}